\documentstyle[11pt]{article}
\def\Bbb R{{\rm \bf R}}
\def\proclaim#1{\vskip2mm{\bf #1}\em}
\def\endproclaim{\em \vskip2mm}
\def\tag#1{\eqno(#1)}
\def\gathered{\begin{array}{c}}
\def\endgathered{\end{array}}
\def\text{\mbox}

\begin{document}

\title {Revisiting the probe and enclosure methods}
\author{Masaru IKEHATA\footnote{
Laboratory of Mathematics,
Graduate School of Advanced Science and Engineering,
Hiroshima University, Higashhiroshima 739-8527, JAPAN}
\footnote{Emeritus Professor at Gunma University, Maebashi 371-8510, JAPAN}}
\maketitle
 
\begin{abstract}
This paper is concerned with reconstruction issue of inverse obstacle problems
governed by partial differential equations
and consists of two parts.

(i) The first part considers the foundation of the probe and enclosure methods for an impenetrable obstacle
embedded in a medium governed by the stationary Schr\"odinger equation.
Under a general framework, some natural estimates for a quantity computed from
a pair of the Dirichlet and Neumann data on the outer surface of the body occupied by the medium are given.
The estimates enables us to derive almost immediately the necessary asymptotic behaviour of indicator functions for both methods.

(ii)  The second one considers the realization of the enclosure method for a penetrable obstacle embedded in an absorbing medium
governed by the stationary Schr\"odinger equation.  
The unknown obstacle considered here is modeled by a perturbation term added to
the background complex-valued $L^{\infty}$-potential.
Under a jump condition on the term across the boundary of the obstacle and some kind of regularity for the obstacle surface including Lipschitz one
as a special case, the enclosure method using infinitely many pairs
of the Dirichlet and Neumann data is established.


\noindent
AMS: 35R30

\end{abstract}


\section{Introduction}

The probe and enclosure methods initiated by the author in \cite{IProbe0}, \cite{IProbe1}, \cite{IProbe2}, \cite{I1}, \cite{IE}, \cite{I2} and \cite{IProbe05}
are methodologies in reconstruction issue for inverse obstacle problems governed by partial differential equations.
Nowadays we have various applications, see the recent survey paper \cite{ISurvey} for the results and references. 
However, looking back to some of those applications, one has not yet fully developed its possibility.

In this paper, we consider the following two topics on inverse obstacle problems governed by
the stationary Schr\"odinger equation:

\noindent
(i)  the foundation of the probe and enclosure methods for an {\it impenetrable obstacle};

\noindent
(ii)  the enclosure method for a {\it penetrable obstacle} embedded 
in an {\it absorbing medium}.

To explain the content of the first part, consider an inverse obstacle problem governed by the Laplace equation in a bounded domain $\Omega$ of $\Bbb R^3$
with smooth boundary.  The problem is motivated by a possibility of application to nondestructive testing
and formulated as follows.
Given an arbitrary solution $v=v(x)$ of the Laplace equation in the whole domain $\Omega$,
let $u=u(x)$ be a solution of 
$$
\left\{
\begin{array}{ll}
\displaystyle
\Delta u=0, & x\in\Omega\setminus\overline{D},\\
\\
\displaystyle
\frac{\partial u}{\partial\nu}=0, & x\in\partial D,
\\
\\
\displaystyle
u=v, & x\in\partial\Omega
\end{array}
\right.
\tag {1.1}
$$
where $D$ is an open subset of $\Bbb R^3$ with, say, smooth boundary and satisfies that $\overline{D}\subset\Omega$
and $\Omega\setminus\overline{D}$ is connected;  $\nu$ denotes the unit outward normal vector to $\partial D$.
The $D$ is a mathematical model of a defect occurred in a body $\Omega$.
Consider the pair of the  Dirichlet and Neumann data $u\vert_{\partial\Omega}=v\vert_{\partial\Omega}$ and
$\frac{\partial u}{\partial\nu}\vert_{\partial\Omega}$, where $\nu$ denotes also the unit outward normal vector to $\partial\Omega$.
The problem is to extract information about the geometry of $D$ from a set of the pairs.

The probe and enclosure methods, say, see section 2 in \cite{ICarleman} and section 5 in \cite{I2}, enables us to extract $D$ itself and the convex hull of $D$, respectively
from the indicator functions calculated from specially chosen infinitely many pairs of the Dirichlet and Neumann data.

The one of common key points for both methods is the fundamental estimates for the pair:
$$\displaystyle
\Vert \nabla v\Vert_{L^2(D)}^2\le \left<\frac{\partial v}{\partial\nu}\vert_{\partial\Omega}-\frac{\partial u}{\partial\nu}\vert_{\partial\Omega}\,,\overline{v}\,\right>
\le C\Vert v\Vert_{H^1(D)}^2,
\tag {1.2}
$$
where $\left<\,,\,\right>$ denotes the dual pairing between $H^{-\frac{1}{2}}(\partial\Omega)$ and $H^{\frac{1}{2}}(\partial\Omega)$
and $C$ a positive constant independent of $v$.

However, from the beginning of those methods, it was not known whether (1.2) type estimates is valid for the Helmholtz equation with a fixed wave number or not.
So to realize the both methods in \cite{IProbe2} and \cite{I2} for inverse obstacle problems governed by the Helmholtz equation
the author developed some technical argument.  The argument worked well for the probe method.
However, for the enclosure method \cite{I2} we needed the argument combined with 
some additional restriction on the curvature of the surface of the obstacle.
Later, in \cite{SY}, \cite{SY2} Sini-Yoshida developed an argument to remove such restriction and obtained the estimates:
$$\displaystyle
C_1\Vert \nabla v\Vert_{L^2(D)}^2-C_2\Vert v\Vert_{L^2(D)}^2\le 
\left<\frac{\partial v}{\partial\nu}\vert_{\partial\Omega}
-\frac{\partial u}{\partial\nu}\vert_{\partial\Omega}\,,\overline{v}\,\right>
\le C_3\Vert v\Vert_{H^1(D)}^2,
\tag {1.3}
$$
where $u$ is a solution of 
$$
\left\{
\begin{array}{ll}
\displaystyle
\Delta u+k^2u=0, & x\in\Omega\setminus\overline{D},\\
\\
\displaystyle
\frac{\partial u}{\partial\nu}=0, & x\in\partial D,
\\
\\
\displaystyle
u=v, & x\in\partial\Omega
\end{array}
\right.
$$
and $v$ is an arbitrary solution of equation $\Delta v+k^2v=0$ in $\Omega$.  It is assumed that $k^2$ is not a Dirichlet eigenvalue of $-\Delta$ in $D$ nor
eigenvalue of $-\Delta$ in $\Omega\setminus\overline{D}$ with the Dirichlet and Neumann boundary conditions on $\partial\Omega$ and $\partial D$, respectively.
Note that the validity of the upper estimate part on (1.3) is easily deduced
and already known at the beginning of the probe and enclosure methods.
Thus, the real problem was to establish the lower estimate part.

The essence of their proof on the lower estimate is as follows.
First we have the integral expression of the middle term on (1.3) from \cite{I2}:
$$\displaystyle
\left<\frac{\partial v}{\partial\nu}\vert_{\partial\Omega}
-\frac{\partial u}{\partial\nu}\vert_{\partial\Omega}\,,\overline{v}\,\right>
=\Vert\nabla v\Vert_{L^2(D)}^2-k^2\Vert v\Vert_{L^2(D)}^2+\Vert\nabla w\Vert_{L^2(\Omega\setminus\overline{D})}^2-k^2\Vert w\Vert_{L^2(\Omega\setminus\overline{D})}^2,
$$
where $w=u-v$.  By dropping the third term we have
$$\displaystyle
\left<\frac{\partial v}{\partial\nu}\vert_{\partial\Omega}
-\frac{\partial u}{\partial\nu}\vert_{\partial\Omega}\,,\overline{v}\,\right>
\ge
\Vert\nabla v\Vert_{L^2(D)}^2-k^2\Vert v\Vert_{L^2(D)}^2-k^2\Vert w\Vert_{L^2(\Omega\setminus\overline{D})}^2.
\tag {1.4}
$$
This is the common part as the previous approach in \cite{I2}
and thus the problem is to give a good estimate for $\Vert w\Vert_{L^2(\Omega\setminus\overline{D})}$
which is {\it weaker} than $\Vert\nabla v\Vert_{L^2(D)}$.
For this the author in \cite{I2} made use of a character of specially chosen input $v$ under some restriction, see Lemma 4.2 in \cite{I2}\footnote{The reason why the author did not pursue deeper at that time was that at that time, the search for a new method was the first, and the pursuit of sharp results was the second.}.
In contrast to this, they deduced that, for a fixed $s\in\,]\frac{3}{2},\,1[$ 
independent of general input $v$
$$
\displaystyle
\Vert w\Vert_{L^2(\Omega\setminus\overline{D})}\le C\Vert v\Vert_{H^s(D)}.
\tag {1.5}
$$
They combined this with the additional estimate (cf. Theorem 1.4.3.3 on page 26 in  \cite{Gr})
$$\displaystyle
\Vert v\Vert_{H^s(D)}^2\le \epsilon\Vert\nabla v\Vert_{L^2(D)}^2+C_{\epsilon}\Vert v\Vert_{L^2(D)}^2,
\tag {1.6}
$$
where $C_{\epsilon}$ is a positive number depending on an arbitrary small positive number $\epsilon$ and independent of $v$.
From (1.4), (1.5) and (1.6) for a sufficiently small $\epsilon$, one gets the lower estimate on (1.3).

The proof of (1.5) which is one of the main parts of their paper employs the potential theoretic approach on fractional order Sobolev spaces and quite involved.
Especially its central part, Lemma 4.3 in \cite{SY} states the invertibility of an operator involving the adjoint of the double layer potential 
on a negative fractional order Sobolev space on $\partial D$.
 
In the first part of this paper we give an alternative and elementary approach, however, which yields the sharper estimate than (1.5), that is,
$$\displaystyle
\Vert w\Vert_{L^2(\Omega\setminus\overline{D})}\le C\Vert v\Vert_{L^2(D)}.
$$
From this we immediately obtain the lower estimate on (1.3)(without using (1.6)).
It is a special case of the result under a {\it general framework}.  The major advantage of our approach is:
we do not need any detailed knowledge of the governing equation for the background medium in contrast to their approach in which the fundamental solution
of the Helmholtz equation played the central role.
Thus it covers inverse problems for impenetrable obstacles governed by the stationary Schr\"odinger equation almost immediately.
As a direct consequence we obtain an extension and improvement of the author's previous application of the enclosure method
to an obstacle with impedance boundary condition in \cite{Iimpedance}.

The description above is the comparison of the methods for the proof of the lower estimate on (1.3).
From a technical point of view, it should be pointed out that,  the proof in \cite{SY,SY2} mentioned above
which is based on the potential theory covers the case when $\partial D$ is Lipschitz. 
Our proof even in this simple situation needs a higher regularity $\partial D\in C^{1,1}$.  This is because of the trace theorem of $H^2$-functions.
See Lemma 2.2.

The content of the second part is as follows.
The unknown obstacle considered here is a penetrable one and modeled by a perturbation term added to
the background {\it complex-valued} $L^{\infty}$-potential of the stationary Schr\"odinger equation.  
Under a condition on the jump of the real or imaginary part of the term across the boundary of the obstacle and some kind of regularity for the obstacle surface
including Lipschitz one
as a special case, the enclosure method using infinitely many pairs of the Dirichlet and Neumann data is established.  This is a full extension of
a result in the unpublished manuscript \cite{IEx} (see also \cite{ISurvey}) in which an inverse obstacle problem
for a non-absorbing medium governed by the stationary Schr\"odiner equation was considered.

Now let us formulate the two problems mentioned above more precisely and describe statements of the results.

\subsection{Impenetrable obstacle}

Let $\Omega$ be a bounded domain of $\Bbb R^3$ with smooth boundary.
Let $v\in H^1(\Omega)$ be a weak solution of
$$\begin{array}{ll}
\displaystyle
\Delta v+V_0(x)v=0, & x\in\Omega,
\end{array}
\tag {1.7}
$$
where $V_0\in L^{\infty}(\Omega)$.
This means that, for all $\psi\in H^1(\Omega)$ with $\psi=0$ on $\partial\Omega$ in the sense of the trace we have
$$\displaystyle
-\int_{\Omega}\nabla v\cdot\nabla\psi\,dx+\int_{\Omega}\,V_0(x)v\psi\,dx=0.
\tag {1.8}
$$

Then the bounded linear functional $\frac{\partial v}{\partial\nu}\vert_{\partial\Omega}\in H^{-\frac{1}{2}}(\partial\Omega)$ is well defined by the formula
$$\displaystyle
\left<\frac{\partial v}{\partial\nu}\vert_{\partial\Omega},\,f
\right>
=\int_{\Omega}\,\nabla v\cdot\nabla\eta\,dx-\int_{\Omega}\,V_0(x)v\eta\,dx,
\tag {1.9}
$$
where $\eta$ is an arbitrary element in $H^1(\Omega)$ such that $\eta=f\in H^{\frac{1}{2}}(\partial\Omega)$.  Note that
the well-definedness means that the value (1.9) does not change when $\eta$ is replaced with any other one having the same trace on $\partial\Omega$.
This is because of (1.8).
The boundedness is a consequence a special choice of $\eta$ such that $\Vert \eta\Vert_{H^1(\Omega)}\le C\Vert f\Vert_{H^{1/2}(\partial\Omega)}$
with a positive constant $C$ independent of $f$, whose existence is ensured by the trace theorem about the lifting (cf., Theorem 1.5.1.3 in \cite{Gr}).  
The functional $\frac{\partial v}{\partial\nu}\vert_{\partial\Omega}$ is boundary local in the sense that
it can be determined by the value of $v$ in $U\cap\Omega$, where $U$ is an arbitrary {\it small} neighbourhood of $\partial\Omega$.
This is also a consequence of (1.8).

Let $D$ be a nonempty open set of $\Bbb R^3$ such that $\overline D\subset\Omega$, $\Omega\setminus\overline D$ is connected and $\partial D$ is Lipschitz.
Let $\lambda\in L^{\infty}(\partial D)$.  The $\lambda$ can be a complex valued function.

Given a weak solution $v$ of (1.7),
let $u\in H^1(\Omega\setminus\overline D)$ be a weak solution of
$$\left\{
\begin{array}{ll}
\displaystyle
\Delta u+V_0(x)u=0, & x\in\partial\Omega\setminus\overline D,\\
\\
\displaystyle
\frac{\partial u}{\partial\nu}+\lambda(x)u=0, & x\in\partial D,
\\
\\
\displaystyle
u=v, & x\in\partial\Omega.
\end{array}
\right.
\tag {1.10}
$$

This means that, for all $\varphi\in H^1(\Omega\setminus\overline D)$ with $\varphi=0$ on $\partial\Omega$ in the sense
of the trace we have
$$\displaystyle
\int_{\partial D}\,\lambda(x)u\varphi\,dS-\int_{\Omega\setminus\overline D}\nabla u\cdot\nabla\varphi\,dx+\int_{\Omega\setminus\overline D}\,V_0(x)u\varphi\,dx=0
\tag {1.11}
$$
and $u=v$ on $\partial\Omega$ in the sense of trace.  Note that $u$ depends on $v$ on $\partial\Omega$.

Then the bounded linear functional $\frac{\partial u}{\partial\nu}\vert_{\partial\Omega}\in H^{-\frac{1}{2}}(\partial\Omega)$ is well defined by the formula
$$\displaystyle
\left<\frac{\partial u}{\partial\nu}\vert_{\partial\Omega},\,f\right>
=-\int_{\partial D}\,\lambda(x)\,u\phi\,dS+\int_{\Omega\setminus\overline D}\,\nabla u\cdot\nabla\phi\,dx-\int_{\Omega\setminus\overline D}\,V_0(x)u\phi\,dx,
\tag {1.12}
$$
where $\phi\in H^1(\Omega\setminus\overline D)$ such that $\phi=f\in H^{\frac{1}{2}}(\partial\Omega)$.
Needless to say, for the well-definedness and boundedness, we have made use of (1.12) and the trace theorem for the lifting \cite{Gr}, respectively.

Our starting point is the following representation formula.

\proclaim{\noindent Proposition 1.1.}  We have
$$\begin{array}{ll}
\displaystyle
\left<\frac{\partial v}{\partial\nu}\vert_{\partial\Omega}-\frac{\partial u}{\partial\nu}\vert_{\partial\Omega},\,\overline{v}\vert_{\partial\Omega}
\right>
&
\displaystyle
=2i\,\int_{\partial D}\,\text{Im}\,(\lambda(x))\,w\,\overline{v}\,dS+2i\int_{\Omega\setminus\overline{D}}\,\text{Im}\,(V_0(x))\,w\overline v\,dx
\\
\\
\displaystyle
&
\displaystyle
\,\,\,
-\int_{\partial D}\,\overline{\lambda(x)}\,\vert w\vert^2\,dS
+\int_{\Omega\setminus\overline D}\,\vert\nabla w\vert^2\,dx
-\int_{\Omega\setminus\overline D}\,\overline{V_0(x)}\vert w\vert^2\,dx
\\
\\
\displaystyle
&
\displaystyle
\,\,\,
+\int_{\partial D}\,\lambda(x)\,\vert v\vert^2\,dS
+\int_D\vert\nabla v\vert^2\,dx
-\int_D\,V_0(x)\vert v\vert^2\,dx,
\end{array}
\tag {1.13}
$$
where $w=u-v$.

\endproclaim
{\it\noindent Proof.}
The $w$ satisfies $w=0$ on $\partial\Omega$ in the sense of the trace.
From (1.9) for $\eta=\overline{v}$ and (1.12) for $\phi=\overline{v}$, respectively we have
$$\begin{array}{l}
\displaystyle
\left<\frac{\partial v}{\partial\nu}\vert_{\partial\Omega},\,\overline{v}\vert_{\partial\Omega}\,\right>
=\int_{\Omega}\nabla v\cdot\nabla\overline{v}\,dx+\int_{\Omega}\,V_0(x)v\overline{v}\,dx.
\end{array}
$$
and
$$\begin{array}{l}
\displaystyle
\left<\frac{\partial u}{\partial\nu}\vert_{\partial\Omega},\,\overline{v}\vert_{\partial\Omega}\,\right>
=-\int_{\partial D}\,\lambda(x)\,(v+w)\overline{v}\,dS+
\int_{\Omega\setminus\overline D}\,\nabla (v+w)\cdot\nabla\overline{v}\,dx-\int_{\Omega\setminus\overline D}\,V_0(x)(v+w)\overline{v}\,dx.
\end{array}
$$
Thus one gets
$$\begin{array}{l}
\displaystyle
\,\,\,\,\,\,
\left<\frac{\partial v}{\partial\nu}\vert_{\partial\Omega}-\frac{\partial u}{\partial\nu}\vert_{\partial\Omega},\,\overline{v}\vert_{\partial\Omega}\right>
\\
\\
\displaystyle
=\int_{\partial D}\,\lambda(x)\,w\overline{v}\,dS
-\int_{\Omega\setminus\overline D}\,\nabla w\cdot\nabla\overline {v}\,dx
+\int_{\Omega\setminus\overline D}\,V_0(x)w\overline{v}\,dx\\
\\
\displaystyle
\,\,\,
+\int_{\partial D}\,\lambda(x)\,\vert v\vert^2\,dS
+\int_D\vert\nabla v\vert^2\,dx
-\int_D\,V_0(x)\vert v\vert^2\,dx
\\
\\
\displaystyle
=\int_{\partial D}\,\lambda(x)\,w\overline{(u-w)}\,dS
-\int_{\Omega\setminus\overline D}\,\nabla w\cdot\nabla\overline {(u-w)}\,dx
+\int_{\Omega\setminus\overline D}\,V_0(x)w\overline{(u-w)}\,dx
\\
\\
\displaystyle
\,\,\,
+\int_{\partial D}\,\lambda(x)\,\vert v\vert^2\,dS
+\int_D\vert\nabla v\vert^2\,dx
-\int_D\,V_0(x)\vert v\vert^2\,dx
\\
\\
\displaystyle
=\int_{\partial D}\,\lambda(x)\,w\overline{u}\,dS
-\int_{\Omega\setminus\overline D}\,\nabla w\cdot\nabla\overline{u}\,dx
+\int_{\Omega\setminus\overline D}\,V_0(x)w\overline{u}\,dx
\\
\\
\displaystyle
\,\,\,
-\int_{\partial D}\,\lambda(x)\,\vert w\vert^2\,dS
+\int_{\Omega\setminus\overline D}\,\vert\nabla w\vert^2\,dx
-\int_{\Omega\setminus\overline D}\,V_0(x)\vert w\vert^2\,dx
\\
\\
\displaystyle
\,\,\,
+\int_{\partial D}\,\lambda(x)\,\vert v\vert^2\,dS
+\int_D\vert\nabla v\vert^2\,dx
-\int_D\,V_0(x)\vert v\vert^2\,dx.
\end{array}
$$
Here using (1.11) for $\varphi=\overline{w}$ we have
$$\begin{array}{l}
\displaystyle
\,\,\,\,\,\,
\int_{\partial D}\,\lambda(x)\,w\overline{u}\,dS
-\int_{\Omega\setminus\overline D}\,\nabla w\cdot\nabla\overline{u}\,dx
+\int_{\Omega\setminus\overline D}\,V_0(x)w\overline{u}\,dx
\\
\\
\displaystyle
=\overline{\int_{\partial D}\,\lambda(x)\,u\overline{w}\,dS
-
\int_{\Omega\setminus\overline D}\,\nabla u\cdot\nabla\overline{w}\,dx
+\int_{\Omega\setminus\overline D}\,V_0(x)u\overline{w}\,dx\,}
\\
\\
\displaystyle
\,\,\,
+\int_{\partial D}(\lambda(x)-\overline{\lambda(x)}\,)\,w\overline{u}\,dS
+\int_{\Omega\setminus\overline{D}}(V_0(x)-\overline{V_0(x)})\,w\overline u\,dx
\\
\\
\displaystyle
=\int_{\partial D}\,(\lambda(x)-\overline{\lambda(x)}\,)\,w\overline{v}\,dS
+\int_{\Omega\setminus\overline{D}}(V_0(x)-\overline{V_0(x)})\,w\overline v\,dx\\
\\
\displaystyle
\,\,\,
+\int_{\partial D}\,(\lambda(x)-\overline{\lambda(x)}\,)\,\vert w\vert^2\,dS
+\int_{\Omega\setminus\overline D}\,(V_0(x)-\overline{V_0(x)}\,\vert w\vert^2\,dx.
\end{array}
$$
Therefore we have obtained the desired formula.

\noindent
$\Box$

Note that Proposition 1.1 is valid under the existence of $u$.  We never make use of any estimate on $u$.

To go further hereafter we impose an assumption concerning with the uniqueness and existence of $u$ together with 
some estimates on the reflected solution $w=u-v$.

$\quad$

{\bf\noindent Assumption 1.}
Given $F\in L^2(\Omega\setminus\overline D)$ there exists a unique weak solution $p\in H^1(\Omega\setminus\overline D)$
of
$$
\left\{
\begin{array}{ll}
\displaystyle
\Delta p+V_0(x)p=F, & x\in\Omega\setminus\overline D,
\\
\\
\displaystyle
\frac{\partial p}{\partial\nu}+\lambda(x)\,p=0, & x\in\partial D,\\
\\
\displaystyle
p=0, & x\in\partial\Omega.
\end{array}
\right.
\tag{1.14}
$$
and that the unique solution satisfies
$$\displaystyle
\Vert p\Vert_{L^2(\Omega\setminus\overline D)}\le C\Vert F\Vert_{L^2(\Omega\setminus\overline D)},
\tag {1.15}
$$
where $C$ is a positive constant independent of $F$.

$\quad$

This means that $p$ satisfies $p=0$ on $\partial\Omega$ in the sense of the trace and,
for all $\phi\in H^1(\Omega\setminus\overline D)$ with $\phi=0$ on $\partial\Omega$ in the sense of the trace, we have
$$
\displaystyle
\int_{\partial D}\,\lambda(x)p\phi\,dS-\int_{\Omega\setminus\overline D}\,\nabla p\cdot\nabla\phi\,dx+\int_{\Omega\setminus\overline D}\,V_0(x)p\phi\,dx=
\int_{\Omega\setminus\overline D}\,F\phi\,dx.
\tag {1.16}
$$
From this with $\phi=\overline{p}$ we have
$$\displaystyle
\Vert\nabla p\Vert_{L^2(\Omega\setminus\overline{D})}^2
=\int_{\partial D}\,\lambda(x)\vert p\vert^2\,dS+\int_{\Omega\setminus\overline{D}}V_0(x)\vert p\vert^2\,dx-\int_{\Omega\setminus\overline{D}}\,F\overline{p}\,dx.
$$
Here using Theorem 1.5.1.10 in \cite{Gr}, we have, for all $\epsilon>0$
$$\begin{array}
{ll}
\displaystyle
\left\vert\int_{\partial D}\,\lambda(x)\vert p\vert^2\,dS\right\vert
&
\displaystyle
\le C\Vert p\vert_{\partial D}\Vert_{L^2(\partial D)}^2
\\
\\
\displaystyle
&
\displaystyle
\le C(\epsilon\Vert\nabla p\Vert_{L^2(\Omega\setminus\overline{D})}^2+\epsilon^{-1}\Vert p\Vert_{L^2(\Omega\setminus\overline{D})}^2\,).
\end{array}
$$
Thus, choosing a small $\epsilon$ and using (1.16), we conclude that 
$$\displaystyle
\Vert\nabla p\Vert_{L^2(\Omega\setminus\overline D)}\le C\Vert F\Vert_{L^2(\Omega\setminus\overline D)}
$$
and (1.15) yields
$$\displaystyle
\Vert p\Vert_{H^1(\Omega\setminus\overline{D})}\le C\Vert F\Vert_{L^2(\Omega\setminus\overline{D})}.
\tag {1.17}
$$

By Assumption 1 together with a standard lifting argument, we know that, given $v$ which is an arbitrary weak solution of (1.7)
there exists a unique weak solution $u\in H^1(\Omega\setminus\overline{D})$ of (1.10) in the sense mentioned above.

$\quad$

{\bf\noindent Problem 1.}  Find the (1.3) like estimates for 
$\text{Re}\,\left<\frac{\partial v}{\partial\nu}\vert_{\partial\Omega}-\frac{\partial u}{\partial\nu}\vert_{\partial\Omega},\,\overline{v}\vert_{\partial\Omega}
\right>$.

$\quad$

Here we impose the following assumption\footnote{It seems that the treatment of the second term of the right-hand side on (1.13) looks difficult in the case when $\text{Im}\,V_0(x)\not=0$.  The reason for the appearance of the term is simple.  The complex conjugate of $v$ satisfies 
$\Delta\overline{v}+\overline{V_0(x)}\overline{v}=0$ not (1.7).
See also the last section of this paper}.

$\quad$

{\bf\noindent Assumption 2.}  $\text{Im}\,V_0(x)=0$ a. e. $x\in\Omega\setminus\overline{D}$.

$\quad$

Then (1.13) becomes
$$\begin{array}{l}
\displaystyle
\,\,\,\,\,\,
\left<\frac{\partial v}{\partial\nu}\vert_{\partial\Omega}-\frac{\partial u}{\partial\nu}\vert_{\partial\Omega},\,\overline{v}\vert_{\partial\Omega}
\right>
\\
\\
\displaystyle
=2i\int_{\partial D}\,\text{Im}\,(\lambda(x))\,w\overline{v}\,dS\\
\\
\displaystyle
\,\,\,
-\int_{\partial D}\,\overline{\lambda(x)}\,\vert w\vert^2\,dS
+
\int_{\Omega\setminus\overline D}\,\vert\nabla w\vert^2\,dx-\int_{\Omega\setminus\overline D}\,V_0(x)\,\vert w\vert^2\,dx
\\
\\
\displaystyle
\,\,\,
+\int_{\partial D}\,\lambda(x)\,\vert v\vert^2\,dS
+\int_D\vert\nabla v\vert^2\,dx-\int_D\,V_0(x)\vert v\vert^2\,dx.
\end{array}
\tag {1.18}
$$
Here we have
$$\begin{array}{ll}
\displaystyle
2i\int_{\partial D}\,\text{Im}\,(\lambda(x))\,w\overline{v}\,dS
&
\displaystyle
=-2\int_{\partial D}\,\text{Im}\,(\lambda(x))\,\text{Im}\,(w\overline{v})\,dS
+2i\int_{\partial D}\,\text{Im}\,(\lambda(x))\,\text{Re}\,(w\overline{v})\,dS.
\end{array}
$$
Thus taking the real part of the both sides on (1.18), we obtain
$$\begin{array}{l}
\displaystyle
\,\,\,\,\,\,
\text{Re}\,\left<\frac{\partial v}{\partial\nu}\vert_{\partial\Omega}-\frac{\partial u}{\partial\nu}\vert_{\partial\Omega},\,\overline{v}\vert_{\partial\Omega}
\right>
\\
\\
\displaystyle
=-2\int_{\partial D}\,\text{Im}\,(\lambda(x))\,\text{Im}\,(w\overline{v})\,dS
\\
\\
\displaystyle
\,\,\,
-\int_{\partial D}\,\text{Re}\,(\lambda(x))\,\vert w\vert^2\,dS+\int_{\Omega\setminus\overline{D}}\,\vert\nabla w\vert^2\,dx-\int_{\Omega\setminus\overline{D}}\,\text{Re}\,(V_0(x))\,\vert w\vert^2\,dx
\\
\\
\displaystyle
\,\,\,
+\int_{\partial D}\,\text{Re}\,(\lambda(x))\,\vert v\vert^2\,dS+\int_D\,\vert\nabla v\vert^2\,dx-\int_D\,\text{Re}\,(V_0(x))\,\vert v\vert^2\,dx.
\end{array}
\tag {1.19}
$$

Now we state one of the main results of this paper.

\proclaim{\noindent Theorem 1.1.}  Assume that $\partial D$ is $C^{1,1}$ and $\lambda\in C^{0,1}(\partial D)$.
Under Assumptions 1 and 2
we have
$$\displaystyle
C_1\Vert\nabla v\Vert_{L^2(D)}^2-C_2\Vert v\Vert_{L^2(D)}^2
\le 
\text{Re}\,
\left<\frac{\partial v}{\partial\nu}\vert_{\partial\Omega}
-\frac{\partial u}{\partial\nu}\vert_{\partial\Omega},\,\overline{v}\vert_{\partial\Omega}
\right>
\le C_3\Vert v\Vert_{H^1(D)}^2,
\tag {1.20}
$$
where $C_1$, $C_2$ and $C_3$ are positive numbers independent of $v$.

\endproclaim

Some remarks are in order.

\noindent
(i)  
The validity of the system of inequalities on (1.20), especially the lower estimate is independent of the signature
of the real or imaginary part of $\lambda$.

\noindent
(ii)  
Taking the imaginary part of the both sides on (1.18), we obtain
$$\begin{array}{ll}
\displaystyle
\text{Im}\,\left<\frac{\partial v}{\partial\nu}\vert_{\partial\Omega}-\frac{\partial u}{\partial\nu}\vert_{\partial\Omega},\,\overline{v}\vert_{\partial\Omega}
\right>
&
\displaystyle
=2\int_{\partial D}\,\text{Im}\,(\lambda(x))\,\text{Re}\,(w\overline{v})\,dS\\
\\
\displaystyle
&
\displaystyle
\,\,\,
+\int_{\partial D}\,\text{Im}\,(\lambda(x))\,\vert w\vert^2\,dS
+\int_{\partial D}\,\text{Im}\,(\lambda(x))\,\vert v\vert^2\,dS
\\
\\
\displaystyle
&
\displaystyle
=\int_{\partial D}\,\text{Im}\,(\lambda(x))\,\vert u\vert^2\,dS.
\end{array}
$$
It would be not suitable for us to extract the information about the location of $D$ itself from this equation
since in Theorem 1.1, there is no restriction on the signature of $\text{Im}\,\lambda$.

\subsection{Penetrable obstacle in an absorbing medium}

Let $\Omega$ be bounded domain of $\Bbb R^3$ with $C^{1,1}$-boundary.
In this subsection we consider an inverse obstacle problem governed by the stationary Schr\"odinger equation
$$\begin{array}{ll}
\displaystyle
\Delta u+V_0(x)u+V(x)u=0,
&
x\in\Omega,
\end{array}
\tag {1.21}
$$
where both $V_0$ and $V$ are {\it complex} valued $L^{\infty}(\Omega)$ functions.

Let $u\in H^1(\Omega)$ be a weak solution of (1.21).
As usual, define the linear functional $\frac{\partial u}{\partial\nu}\vert_{\partial\Omega}$ on $H^{\frac{1}{2}}(\partial\Omega)$ by the formula
$$\displaystyle
\left<\frac{\partial u}{\partial\nu}\vert_{\partial\Omega},f\right>
=\int_{\Omega}\,\nabla u\cdot\nabla \eta\,dx-\int_{\Omega}\,(V_0(x)+V(x))\,u\eta\,dx,
$$
where $\eta$ is an arbitrary element in $H^1(\Omega)$ such that $\eta=f\in H^{\frac{1}{2}}(\partial\Omega)$
and $\nu$ denotes the unit outward normal vector field to $\partial\Omega$.
Note that this right-hand side is invariant to any $\eta$ as long as the condition $\eta=f$ on $\partial\Omega$ is satisfied.
Thus choosing a special lifting for $f$, we see that this functional is bounded, that is, belongs to $H^{-\frac{1}{2}}(\partial\Omega)$.

We consider the following problem.

$\quad$

{\bf\noindent Problem 2.} 
Assume that $V_0$ is known.  Extract information about the {\it discontinuity}, roughly speaking the place where $V(x)\not=0$ from 
the Neumann data $\frac{\partial u}{\partial\nu}\vert_{\partial\Omega}$
for the weak solution of (1.21) with some known Dirichlet data $u\vert_{\partial\Omega}$.

$\quad$

Our aim is not to reconstruct the full knowledge of $V$ or $V_0+V$ itself unlike \cite{SU2} and \cite{Na}.
Here we employ the enclosure method introduced in \cite{IE}\footnote{For an application of the probe method to the case when $V_0(x)\equiv 0$ a.e. $x\in\Omega$,
see Sections 4 and 5 in \cite{IProbe0}.}.

The following assumption for $V_0+V$ corresponds to Assumption 1.

$\quad$

\noindent
{\bf Assumption 3.}  Given an arbitrary {\it complex valued} function $F\in L^2(\Omega)$ there exists a unique weak solution $p\in H^1(\Omega)$ of
$$\left\{
\begin{array}{ll}
\Delta p+V_0(x)p+V(x)p=F, & x\in\Omega,\\
\\
\displaystyle
p(x)=0, & x\in\partial\Omega
\end{array}
\right.
\tag {1.22}
$$
and that the unique solution satisfies
$$\displaystyle
\Vert p\Vert_{L^2(\Omega)}\le C\Vert F\Vert_{L^2(\Omega)},
\tag {1.23}
$$
where $C$ is a positive constant independent of $F$.

$\quad$

As a simple consequence of the weak formulation of (1.22), the estimate (1.23) is equivalent to the following estimate:
$$\displaystyle
\Vert p\Vert_{H^1(\Omega)}\le C'\Vert F\Vert_{L^2(\Omega)},
\tag {1.24}
$$
where the positive constant $C'$ is independent of $F$.
Thus Assumption 3 means that the homogeneous Dirichlet problem (1.22) is well-posed.
Then, by using a lifting argument and (1.24), we see that, given an arbitrary complex valued function $g\in H^{1/2}(\partial\Omega)$ there exists a unique weak 
solution $u\in H^1(\Omega)$ of equation (1.21) with $u=g$ on $\partial\Omega$ in the trace sense.

Succeeding to Assumption 3, to clarify the meaning of discontinuity in Problem 2 we introduce the following assumption.  

$\quad$

\noindent
{\bf Assumption 4.}  There exists a Lebesgue measurable set $D$ of $\Bbb R^3$ having a positive measure and satisfy $D\subset\Omega$
such that the $V(x)$ vanishes for almost all $x\in\Omega\setminus D$.

$\quad$

In contrast to impenetrable case we do not impose Assumption 2.
The set $D$ is a candidate of the discontinuity to be detected.
Note that the condition $D\subset\Omega$ contains the case $\partial D\cap\partial\Omega\not=\emptyset$.
So it is different from the usual condition $\overline D\subset\Omega$.

First we show that under Assumptions 3-4 and $\text{Re}\,V$ has a kind of jump condition across $\partial D$ specified later,
the information about the convex hull of $D$ 
can be obtained by knowing $\frac{\partial u}{\partial\nu}\vert_{\partial\Omega}$
for $g=v\vert_{\partial\Omega}$, where $v$ is an explicit solution of equation (1.7) which
is just the so called complex geometrical optics solution constructed in \cite{SU2} to give an answer 
to the uniqueness issue of the Calder\'on problem \cite{Cal} given by the following steps.

$\quad$

\noindent
1.  Given real unit vector $\omega$ choose $\vartheta$ an arbitrary real unit vector perpendicular to $\omega$.

\noindent
2.  Let $z=\tau(\omega+i\vartheta)$ with $\tau>0$\footnote{Exactly speaking, this choice of $z$ with $z\cdot z=0$ is different from that of \cite{SU2} and simple.}.

\noindent
3.  Compute the unique solution $\Psi(x,z)$ in the weighted 
$L^2$-space $L^2_{-(1-\delta)}(\Bbb R^3)$ for a fixed $\delta\in\,]0,\,1[$ of the integral equation for a sufficiently large $\tau$:
$$\displaystyle
\Psi(x,z)=-\int_{\Bbb R^3}G_z(x-y)\tilde{V_0}(y)\Psi(y,z)\,dy-\int_{\Bbb R^3}G_z(x-y)\tilde{V_0}(y)\,dy,
\tag {1.25}
$$
where $G_z(x)$ is the tempered distribution on $\Bbb R^3$ defined by
$$\displaystyle
G_z(x)=\frac{1}{(2\pi)^3}
\int_{\Bbb R^3}\frac{e^{ix\cdot\xi}}{\vert\xi\vert^2-2iz\cdot\xi}\,d\xi
$$
and $\tilde{V_0}$ denotes the zero extension of $V_0$ outside $\Omega$.
The distribution $e^{x\cdot z}G_z(x)$ on $\Bbb R^3$ is called the Faddeev Green's function at $z\cdot z=0$.

\noindent
4. Define
$$\begin{array}{lll}
\displaystyle
v(x,z)=e^{x\cdot z}(1+\Psi(x,z)), & x\in\Bbb R^3, & \tau>>1.
\end{array}
\tag {1.26}
$$

$\quad$

By the Sylvester-Uhlmann-Nachman estimates \cite{SU2}, \cite{Na} which are concerned with the convolution operator
$f\longmapsto G_z*f$, we have the unique solvability of
the integral equation (1.25) in $L^2_{-(1-\delta)}(\Bbb R^3)$ for $\tau>>1$, the local regularity $\Psi(\,\cdot\,,z)\vert_{\Omega}\in H^2(\Omega)$ and
$\Vert \Psi(\,\cdot\,,z)\vert_{\Omega}\Vert_{H^j(\Omega)}=O(\tau^{j-1})$, for $j=0,1,2$.
In particular, the $v(\,\cdot\,,z)$ given by (1.26) satisfies $v(\,\cdot\,,z)\vert_{\Omega}\in H^2(\Omega)$ and (1.7).
Besides, by the Ramm estimate \cite{R} (14) on page 570 which is concerned with a local estimate of
$(G_z*f)\vert_{\Omega}$, we have\footnote{More precisely he proved the estimate
$$
\displaystyle
\Vert\Psi(\,\cdot\,,z)\vert_{\Omega}\Vert_{L^{\infty}(\Omega)}=O(\frac{\log\tau}{\sqrt{\tau}}).
$$
The decaying order is not important for our purpose.}
$\lim_{\tau\rightarrow\infty}\Vert\Psi(\,\cdot\,,z)\vert_{\Omega}\Vert_{L^{\infty}(\Omega)}=0$.
Thus there exists a positive constant $C$ such that, for all $\tau>>1$ we have
$$\begin{array}{ll}
\displaystyle
C_1e^{\tau x\cdot\omega}\le \vert v(x,z)\vert\le C_2 e^{\tau x\cdot\omega}, & \text{a.e.} x\in\Omega.
\end{array}
\tag {1.27}.
$$
This {\it pointwise-estimate} is the key for our argument since no regularity condition about $V_0$ and $V$ more than their essential boundedness is imposed.

Now having $v=v(x,z)$ given by (1.26), for all $\tau>>1$ we introduce the indicator function in the enclosure method
$$\displaystyle
I_{\omega,\vartheta}\,(\tau)
=\left<\frac{\partial v}{\partial\nu}\vert_{\partial\Omega}-\frac{\partial u}{\partial\nu}\vert_{\partial\Omega}\,, v^*\vert_{\partial\Omega}\,\right>,
\tag {1.28}
$$
where $u$ is the weak solution of (1.21) with $u=v$ on $\partial\Omega$ and 
$$\displaystyle
v^*(x)=v(x,\overline{z}).
\tag {1.29}
$$
Note that, in general $v^*\not=\overline{v}$ and still $v^*$ satisfies (1.7).  
This is the key point to treat complex-valued potential $V_0$.

$\quad$

\noindent
{\bf\noindent Definition 1.1.}
We say that the a real-valued function $K(x)$, $x\in\Omega$ has a positive/negative jump on $\partial D$ from the direction $\omega$ if
there exist $C=C(\omega)>0$ and $\delta=\delta(\omega)>0$ such that $K(x)\ge C$ for almost all $x\in D_{\omega}(\delta)$/
$-K(x)\ge C$ for almost all $x\in D_{\omega}(\delta)$.  Here 
$D_{\omega}(\delta)=D\cap\{x\cdot\omega>h_D(\omega)-\delta\}$ and $h_D(\omega)=\sup_{x\in D}\,x\cdot\omega$ is called the support function of $D$.

$\quad$

Let $\omega$ be a unit vector.   Then, for almost all $s\in\Bbb R$, the set
$$\displaystyle
S_{\omega}(s)=D\cap\{x\in\Bbb R^3\,\vert x\cdot\omega=h_D(\omega)-s\}
$$
is Lebesgue measurable with respect to the two dimensional Lebesgue measure $\mu_2$ on the plane
$x\cdot\omega=h_D(\omega)-s$ identified with the two dimensional Euclidean space $\Bbb R^2$.

$\quad$

{\bf\noindent Definition 1.2.}  Let $1\le p<\infty$.  We say that $D$ is $p$-regular with respect to $\omega$ if there exist
positive numbers $C$ and $\delta$ such that the set $S_{\omega}(s)$
satisfies
$$\begin{array}{ll}
\displaystyle
\mu_2(S_{\omega}(s))\ge Cs^{p-1}, & \text{a.e.}\,s\in]0,\,\delta[.
\end{array}
$$

$\quad$

This type of notion has been used in \cite{IE} for the original version of the enclosure method for the equation
$\nabla\cdot\gamma\nabla u=0$.

Our first result of this subsection under Assumptions 3-4 is the following.

\proclaim{\noindent Theorem 1.2.}  Let $1\le p<\infty$ and $\omega$ be a unit vector.
Assume that $D$ is $p$-regular with respect to $\omega$.
If $\text{Re}\,V$ has a positive or negative jump on $\partial D$ from direction $\omega$, then we have,
for all $\gamma\in\Bbb R$
$$\displaystyle
\lim_{\tau\rightarrow\infty}\,e^{-2\tau t}\tau^{\gamma}\text{Re}\,I_{\omega,\vartheta}(\tau)
=
\left\{
\begin{array}{ll}
\displaystyle
0 & \text{if $t>h_D(\omega)$,}
\\
\\
\displaystyle
\infty & \text{if $t<h_D(\omega)$ and $\text{Re}\,V$ has a positive jump,}
\\
\\
\displaystyle
-\infty & \text{if $t<h_D(\omega)$ and $\text{Re}\,V$ has a negative jump}
\end{array}
\right.
$$
and
$$\begin{array}{ll}
\displaystyle
\liminf_{\tau\rightarrow\infty}e^{-2\tau t}\tau^p\vert\text{Re}\,I_{\omega,\vartheta}(\tau)\vert>0 & \text{if $t=h_D(\omega)$.}
\end{array}
$$
Besides, we have the one-line formula
$$\displaystyle
\lim_{\tau\rightarrow\infty}\,\frac{1}{2\tau}\,\log\vert \text{Re}\,I_{\omega,\vartheta}(\tau)\vert=h_D(\omega).
$$

\endproclaim

\noindent
Note that, as a direct consequence we have the classical characterization 
$$\displaystyle
h_D(\omega)=\inf\{t\in\Bbb R\,\vert\, \lim_{\tau\rightarrow\infty}\,e^{-2\tau t}\text{Re}\,I_{\omega,\vartheta}(\tau)=0\,\}.
$$

In \cite{IE} it is pointed out that if $D$ is a nonempty open set of $\Bbb R^3$ such that $\partial D$ is Lipschitz
or satisfies the interior cone condition, then
$D$ is $3$-regular with respect to an arbitrary unit vector $\omega$.  Thus one gets the following corollary.

\proclaim{\noindent Corollary 1.1.}
Let $D$ in Assumption 3 be an open set of $\Bbb R^3$ such that $D\subset\Omega$ and $\partial D$ is Lipschitz.
If $\text{Re}\,V$ has a positive or negative jump on $\partial D$ from direction $\omega$, 
we have the same conclusion as Theorem 1.2 with $p=3$.

\endproclaim

Concerning with the extraction of discontinuity of $\text{Im}\,V$ which is the main subject of this paper, under Assumptions 3 and 4 we obtain the following result.

\proclaim{\noindent Theorem 1.3.} Let $1\le p<\infty$ and $\omega$ be a unit vector.
Assume that $D$ is $p$-regular with respect to $\omega$.
If $\text{Im}\,V$ has a positive or negative jump on $\partial D$ from direction $\omega$, then we have, for all $\gamma\in\Bbb R$
$$\displaystyle
\lim_{\tau\rightarrow\infty}\,e^{-2\tau t}\tau^{\gamma}\text{Im}\,I_{\omega,\vartheta}(\tau)
=
\left\{
\begin{array}{ll}
\displaystyle
0 & \text{if $t>h_D(\omega)$,}
\\
\\
\displaystyle
\infty & \text{if $t<h_D(\omega)$ and $\text{Im}\,V$ has a positive jump,}
\\
\\
\displaystyle
-\infty & \text{if $t<h_D(\omega)$ and $\text{Im}\,V$ has a negative jump}
\end{array}
\right.
$$
and
$$\begin{array}{ll}
\displaystyle
\liminf_{\tau\rightarrow\infty}e^{-2\tau t}\tau^p\vert\text{Im}\,I_{\omega,\vartheta}(\tau)\vert>0 & \text{if $t=h_D(\omega)$.}
\end{array}
$$
Besides, we have the one-line formula
$$\displaystyle
\lim_{\tau\rightarrow\infty}\,\frac{1}{2\tau}\,\log\vert \text{Im}\,I_{\omega,\vartheta}(\tau)\vert=h_D(\omega).
$$

\endproclaim

\proclaim{\noindent Corollary 1.2.}
Let $D$ in Assumption 3 be an open set of $\Bbb R^3$ such that $D\subset\Omega$ and $\partial D$ is Lipschitz.
If $\text{Im}\,V$ has a positive or negative jump on $\partial D$ from direction $\omega$, 
we have the same conclusion as Theorem 1.3 with $p=3$.

\endproclaim

This paper is organized as follows.
In Section 2 the proof of Theorem 1.1 is given.   It is based on (1.19) and Lemmas 2.1-2.2 in which the estimetes 
of $H^1(\Omega\setminus\overline{D})$ and $L^2(\Omega\setminus\overline{D})$-norms of the reflected solution $w=u-v$ are
given.  The most emphasized one is Lemma 2.2, which ensures $\Vert w\Vert_{L^2(\Omega\setminus\overline{D})}\le C\Vert v\Vert_{L^2(D)}$
for a positive constant independent of $v$.  This estimate has not been known in the previous studies in the probe and enclosure methods.
Section 3 is devoted to the proof of Theorems 1.2-1.3.  It is based on the Alessandrini identity and Lemmas 3.1-3.2 in which 
an upper bound of the $L^2(\Omega)$-norm of the reflected solution $w=u-v$ and the comparison of $\Vert e^{\tau x\cdot\omega}\Vert_{L^1(D)}$
relative to $\Vert e^{\tau x\cdot\omega}\Vert_{L^2(D)}$ as $\tau\rightarrow\infty$, respectively are given.
In particular, the proof of Lemma 3.2 is simple, however, covers a general obstacle compared with the previous unpublished work \cite{IEx} 
and no restriction on the direction $\omega$.
Concerning with the proof of Lemma 3.2, in Section 4 an approach in \cite{SY} for Lipschitz obstacle case is presented.
It is pointed out that their approach does not work in the Lipschitz obstacle case and needs higher regularity.
Section 5 is devoted to applications of Theorem 1.1 to the probe and enclosure methods and an important example covered by Theorems 1.2 and 1.3.
In Section 6, concerning with Assumption 2 and Sini-Yoshida's another result in \cite{SY}, we describe two open problems on the enclosure method.

\section{Proof of Theorem 1.1}

\subsection{Two lemmas}

In this subsections we give two estimates for $w=u-v$ which are crucial for establishing the fundamental inequalities.

\proclaim{\noindent Lemma 2.1.}
Assume that $\partial D$ is Lipschitz and $\lambda\in L^{\infty}(\partial D)$.
Under Assumption 1 we have
$$\displaystyle
\Vert u-v\Vert_{H^1(\Omega\setminus\overline{D})}\le C\Vert v\Vert_{H^1(D)}.
$$

\endproclaim
{\it\noindent Proof.}
The $w$ satisfies $w=0$ on $\partial\Omega$ in the sense of the trace.
Let $F=\overline{w}$ and $p$ be the weak solution in Assumption 2.
Since $u=w+v$ satisfies (1.11) for all $\varphi\in H^1(\Omega\setminus\overline{D})$ with $\varphi=0$ on $\partial D$,
one can substitute $\varphi=p$ into it, we obtain
$$\begin{array}{l}
\displaystyle
\,\,\,\,\,\,
\int_{\partial D}\,\lambda(x)w p\,dS-\int_{\Omega\setminus\overline D}\nabla w\cdot\nabla p\,dx+\int_{\Omega\setminus\overline D}\,V_0(x)w p\,dx
\\
\\
\displaystyle
=-\int_{\partial D}\,\lambda(x)v p\,dS+\int_{\Omega\setminus\overline D}\nabla v\cdot\nabla p\,dx-\int_{\Omega\setminus\overline D}\,V_0(x)v p\,dx.
\end{array}
$$
On the other hand $p$ satisfies (1.16) with $F=\overline{w}$  for $\phi=w$, we obtain
$$\displaystyle
\int_{\partial D}\,\lambda(x)pw\,dS-\int_{\Omega\setminus\overline{D}}\nabla p\cdot\nabla w\,dx+\int_{\Omega\setminus\overline{D}}\,V_0(x)pw\,dx=
\int_{\Omega\setminus\overline{D}}\,\vert w\vert^2\,dx.
$$
Thus we have the expression
$$\begin{array}{ll}
\displaystyle
\Vert w\Vert_{L^2(\Omega\setminus\overline{D})}^2
&
\displaystyle
=-\int_{\partial D}\,\lambda(x)v p\,dS+\int_{\Omega\setminus\overline D}\nabla v\cdot\nabla p\,dx-\int_{\Omega\setminus\overline D}\,V_0(x)v p\,dx.
\end{array}
\tag {2.1}
$$

Here we claim that
$$\displaystyle
\left\vert\int_{\Omega\setminus\overline{D}}\,\nabla v\cdot\nabla p\,dx-\int_{\Omega\setminus\overline{D}}\,V_0(x)vp\,dx\right\vert
\le C(1+\Vert V_0\Vert_{L^{\infty}(D)})\,\Vert v\Vert_{H^1(D)}\Vert p\vert_{\partial D}\Vert_{H^\frac{1}{2}(\partial D)}.
\tag {2.2}
$$
This is proved as follows.
By the lifting in the trace theorem (Theorem 1.5.1.3 in \cite{Gr}), one can find a $\tilde{p}\in H^1(D)$ such that $p=\tilde{p}$ on $\partial D$
in the sense of the trace and 
$$\displaystyle
\Vert\tilde{p}\Vert_{H^1(D)}\le C(D)\Vert p\vert_{\partial D}\Vert_{H^{\frac{1}{2}}(\partial D)}.
\tag {2.3}
$$
Define 
$$\displaystyle
Z=
\left\{
\begin{array}{ll}
\displaystyle
p, & x\in\Omega\setminus\overline{D},
\\
\\
\displaystyle
\tilde{p}, & x\in D.
\end{array}
\right.
$$
By the trace theorem, we have $Z\in H^1(\Omega)$ and $Z=0$ on $\partial\Omega$ in the sense of the trace.

Substituting $\psi=Z$ into (1.8), we obtain
$$\displaystyle
\int_{\Omega\setminus\overline{D}}\,\nabla v\cdot\nabla p\,dx-\int_{\Omega\setminus\overline{D}}\,V_0(x)v p\,dx
=-\int_{D}\,\nabla v\cdot\nabla\tilde{p}\,dx+\int_{D}\,V_0(x)v\tilde{p}\,dx.
\tag {2.4}
$$
A combination of (2.3) and (2.4) yields (2.2).

Now from (2.1) and (2.2) and the trace theorem we obtain
$$\begin{array}{ll}
\displaystyle
\Vert w\Vert_{L^2(\Omega\setminus\overline{D})}^2
&
\displaystyle
\le C(\Vert v\vert_{\partial D}\Vert_{L^2(\partial D)}\Vert p\vert_{\partial D}\Vert_{L^2(\partial D)}
+\Vert v\Vert_{H^1(D)}\Vert p\vert_{\partial D}\Vert_{H^{\frac{1}{2}}(\partial D)}\,)
\\
\\
\displaystyle
&
\displaystyle
\le C\Vert v\Vert_{H^1(D)}\Vert p\Vert_{H^1(\Omega\setminus\overline{D})}.
\end{array}
$$
Then (1.17) with $F=\overline{w}$ yields the rough estimate
$$\displaystyle
\Vert w\Vert_{L^2(\Omega\setminus\overline{D})}\le C\Vert v\Vert_{H^1(D)}.
\tag {2.5}
$$

Next substituting $\varphi=\overline{w}$ into (1.11), we obtain
$$\begin{array}{l}
\displaystyle
\,\,\,\,\,\,
\int_{\partial D}\,\lambda(x)\vert w\vert^2\,dS-\int_{\Omega\setminus\overline D}\vert\nabla w\vert^2\,dx+\int_{\Omega\setminus\overline D}\,V_0(x)\vert w\vert^2\,dx
\\
\\
\displaystyle
=-\int_{\partial D}\,\lambda(x)v\overline{w}\,dS+\int_{\Omega\setminus\overline D}\nabla v\cdot\nabla\overline{w}\,dx-\int_{\Omega\setminus\overline D}\,V_0(x)v\overline{w}\,dx,
\end{array}
$$
that is
$$\begin{array}{ll}
\displaystyle
\Vert\nabla w\Vert_{L^2(\Omega\setminus\overline{D})}^2
&
\displaystyle
=\int_{\partial D}\,\lambda(x)\vert w\vert^2\,dS+\int_{\Omega\setminus\overline D}\,V_0(x)\vert w\vert^2\,dx
\\
\\
\displaystyle
&
\displaystyle
\,\,\,
+\int_{\partial D}\,\lambda(x)v\overline{w}\,dS-\int_{\Omega\setminus\overline D}\nabla v\cdot\nabla\overline{w}\,dx+\int_{\Omega\setminus\overline D}\,V_0(x)v\overline{w}\,dx.
\end{array}
\tag {2.6}
$$
Using similar argument for the derivation of (2.2), we have
$$\displaystyle
\left\vert
\int_{\Omega\setminus\overline D}\nabla v\cdot\nabla\overline{w}\,dx-\int_{\Omega\setminus\overline D}\,V_0(x)v\overline{w}\,dx
\right\vert
\le C\Vert v\Vert_{H^1(D)}\,\Vert w\vert_{\partial D}\Vert _{H^\frac{1}{2}(\partial D)}.
$$
Thus this, (2.6), the trace theorem and Theorem 1.5.1.10 on page 41 in \cite{Gr} yield
$$\begin{array}{ll}
\Vert\nabla w\Vert_{L^2(\Omega\setminus\overline{D})}^2
&
\displaystyle
\le C(\Vert w\vert_{\partial_D}\Vert_{L^2(\partial D)}^2+\Vert w\Vert_{L^2(\Omega\setminus\overline{D})}^2
+\Vert v\Vert_{H^1(D)}\Vert w\Vert_{H^1(\Omega\setminus\overline D)}\,)\\
\\
\displaystyle
&
\displaystyle
\le C(\epsilon\Vert\nabla w\Vert_{L^2(\Omega\setminus\overline{D})}^2+
C_{\epsilon}\Vert w\Vert_{L^2(\Omega\setminus\overline{D})}^2+\epsilon\Vert w\Vert_{H^1(\Omega\setminus\overline{D})}^2+\epsilon^{-1}\Vert\nabla v\Vert_{H^1(D)}^2\,).
\end{array}
$$
Now choosing $\epsilon$ so small, we obtain
$$\displaystyle
\Vert\nabla w\Vert_{L^2(\Omega\setminus\overline{D})}^2\le C(\Vert w\Vert_{L^2(\Omega\setminus\overline{D})}^2+\Vert v\Vert_{H^1(D)}^2\,).
$$
This together with (2.5) yields the desired estimate.

\noindent
$\Box$

We can say more about the regularity of $p$ if $\partial D$ is $C^{1,1}$.   Rewrite (1.14) as
$$
\left\{
\begin{array}{ll}
\displaystyle
\Delta p-p=\tilde{F}, & x\in\Omega\setminus\overline D,
\\
\\
\displaystyle
\frac{\partial p}{\partial\nu}=g, & x\in\partial D,\\
\\
\displaystyle
p=0, & x\in\partial\Omega,
\end{array}
\right.
\tag{2.7}
$$
where
$$\left\{
\begin{array}{ll}
\displaystyle
\tilde{F}=F-(V_0(x)+1)p, & x\in\Omega\setminus\overline{D},
\\
\\
\displaystyle g=-\lambda(x)p, & x\in\partial D.
\end{array}
\right.
$$
Assume that $\lambda\in C^{0,1}(\partial D)$.  By Theorem 6.2.4 on p.277 and Theorem 6.2.6 on p.278 in \cite{Gr} combined with a cut-off function, there exists a 
$\tilde{\lambda}\in C^{0,1}(\Bbb R^3)$ with compact support such that $\lambda=\tilde{\lambda}$ on $\partial D$.
Then, by Theorem 1.4.1.1 on p.21 in \cite{Gr}, we have $\tilde{\lambda} p\in H^1(\Omega\setminus\overline{D})$
and $\Vert\tilde{\lambda} p\Vert_{H^1(\Omega\setminus\overline{D})}\le C\Vert p\Vert_{H^1(\Omega\setminus\overline{D})}$.
Thus this yields $g\in H^{\frac{1}{2}}(\partial D)$ and
$$\displaystyle
\Vert g\Vert_{H^{\frac{1}{2}}(\partial D)}
\le C\Vert p\Vert_{H^1(\Omega\setminus\overline{D})}.
\tag {2.8}
$$

Choose a function $\eta\in C_0^{\infty}(\Bbb R^3)$ in such a way that $\eta\equiv 1$ and $\eta\equiv 0$ in neighbourhoods of
$\partial D$ and $\partial\Omega$, respectively.  Define $p'=\eta p\in H^1(\Omega\setminus\overline{D})$.  The $p'$ is a weak solution of
$$
\left\{
\begin{array}{ll}
\displaystyle
\Delta z-z=G, & x\in\Omega\setminus\overline D,
\\
\\
\displaystyle
\frac{\partial z}{\partial\nu}=g, & x\in\partial D,\\
\\
\displaystyle
\frac{\partial z}{\partial\nu}=0, & x\in\partial\Omega,
\end{array}
\right.
\tag{2.9}
$$
where $G=\eta\tilde{F}+(\Delta\eta)p+2\nabla\eta\cdot\nabla p$.

Corollary 2.2.2.6 on p.92 in \cite{Gr} ensures there exists the unique weak solution $p''$ of (2.9) which belongs to $H^2(\Omega\setminus\overline{D})$
and satisfies
$$\displaystyle
\Vert p''\Vert_{H^2(\Omega\setminus\overline{D})}\le C(\Vert G\Vert_{L^2(\Omega\setminus\overline{D})}+\Vert g\Vert_{H^{\frac{1}{2}}(\partial D)}\,).
$$
Thus $p'$ coincides with $p''$ almost everywhere in $\Omega\setminus\overline{D}$ and
satisfies
$$\begin{array}{ll}
\displaystyle
\Vert p'\Vert_{H^2(\Omega\setminus\overline{D})}
&
\displaystyle
\le C(\Vert\eta\tilde{F}+(\Delta\eta)p+2\nabla\eta\cdot\nabla p\Vert_{L^2(\Omega\setminus\overline{D})}+\Vert g\Vert_{H^{\frac{1}{2}}(\partial D)}\,)
\\
\\
\displaystyle
&
\displaystyle
\le
C(\Vert F\Vert_{L^2(\Omega\setminus\overline{D})}+\Vert p\Vert_{H^1(\Omega\setminus\overline{D})}\,).
\end{array}
$$
Then (1.15) yields
$$\displaystyle
\Vert \eta p\Vert_{H^{2}(\Omega\setminus\overline{D})}
\le C\Vert F\Vert_{L^2(\Omega\setminus\overline{D})}.
\tag {2.10}
$$

Now we are ready to prove a more accurate estimate than that of Lemma 2.1.

\proclaim{\noindent Lemma 2.2.}
Assume that $\partial D$ is $C^{1,1}$ and $\lambda\in C^{0,1}(\partial D)$.
Under Assumption 1 we have
$$\displaystyle
\Vert u-v\Vert_{L^2(\Omega\setminus\overline{D})}\le C\Vert v\Vert_{L^2(D)}.
$$

\endproclaim
{\it\noindent Proof.}
Let $F=\overline{w}$ and $p$ be the weak solution in Assumption 1 again.
It follows from (2.10) with $F=\overline{w}$ we obtain
$$\displaystyle
\Vert \eta p\Vert_{H^{2}(\Omega\setminus\overline{D})}\le C\Vert w\Vert_{L^2(\Omega\setminus\overline D)}.
\tag {2.11}
$$
Next, by the trace theorem on the $C^{1,1}$ domain ( c.f., Theorem 1.5.1.2 in \cite{Gr} ) and (2.8)
one can choose $\tilde{p}\in H^2(D)$ such that
$$\begin{array}{lll}
\displaystyle
\tilde{p}=p,
&
\displaystyle
\frac{\partial\tilde{p}}{\partial\nu}=-\lambda(x)\,p,
&
\displaystyle
x\in\partial D,
\end{array}
$$
in the sense of the trace and satisfies
$$\begin{array}{ll}
\displaystyle
\Vert \tilde{p}\Vert_{H^2(D)}
&
\displaystyle
\le C(\Vert p\vert_{\partial D}\Vert_{H^{\frac{3}{2}}(\partial D)}+\Vert \lambda(x)p\vert_{\partial D}\Vert_{H^{\frac{1}{2}}(\partial D)}\,)\\
\\
\displaystyle
&
\displaystyle
\le C\Vert\eta p\Vert_{H^2(\Omega\setminus\overline{D})}.
\end{array}
$$
This together with (2.11) gives
$$\displaystyle
\Vert \tilde{p}\Vert_{H^2(D)}
\le C\Vert w\Vert_{L^2(\Omega\setminus\overline D)}.
\tag {2.12}
$$
Define
$$\displaystyle
p^*(x)
=
\left\{
\begin{array}{ll}
\displaystyle
p(x), & x\in\,\Omega\setminus\overline D,\\
\\
\displaystyle
\tilde{p}(x), & x\in\,D.
\end{array}
\right.
$$
Then, by the trace theorem, we have $p^*\in H^1(\Omega)$ and $p^*=0$ on $\partial\Omega$.

Substituting $\psi=p^*$ into (1.8), we obtain
$$\displaystyle
-\int_{\Omega\setminus\overline D}\nabla v\cdot\nabla p\,dx+\int_{\Omega\setminus\overline D}\,V_0(x)vp\,dx
=
\int_{D}\nabla v\cdot\nabla\tilde{p}\,dx-\int_{D}\,V_0(x)v\tilde{p}\,dx.
$$
Here, by the interior elliptic regularity of $v$ in $\Omega$ we have $v\vert_{D}\in H^2(D)$.
Then integration by parts \cite{Gr} yields
$$\begin{array}{ll}
\displaystyle
\int_{D}\nabla v\cdot\nabla\tilde{p}\,dx-\int_{D}\,V_0(x)v\tilde{p}\,dx
&
\displaystyle
=\int_{\partial D}\frac{\partial v}{\partial\nu}\,\tilde{p}\,dx-\int_{D}(\Delta v+V_0(x)v)\,\tilde{p}\,dx
\\
\\
\displaystyle
&
\displaystyle
=\int_{\partial D}\frac{\partial v}{\partial\nu}\,\tilde{p}\,dx.
\end{array}
$$
Thus (2.1) becomes
$$\begin{array}{ll}
\displaystyle
\Vert w\Vert_{L^2(\Omega\setminus\overline{D})}^2
&
\displaystyle
=-\int_{\partial D}\,\lambda(x)v p\,dS+\int_{\Omega\setminus\overline D}\nabla v\cdot\nabla p\,dx-\int_{\Omega\setminus\overline D}\,V_0(x)v p\,dx
\\
\\
\displaystyle
&
\displaystyle
=-\int_{\partial D}\,\left(\frac{\partial v}{\partial\nu}+\lambda(x)v\,\right)\tilde{p}\,dS.
\end{array}
\tag {2.13}
$$
Recalling $\tilde{p}\in H^2(D)$, we have
$$\begin{array}{ll}
\displaystyle
\int_{\partial D}\frac{\partial v}{\partial\nu}\,\tilde{p}\,dS
&
\displaystyle
=\int_D\Delta v\tilde{p}\,dx+\int_D\nabla v\cdot\nabla\tilde{p}\,dx\\
\\
\displaystyle
&
\displaystyle
=-\int_D V_0(x)v \tilde{p}\,dx+\int_{\partial D}\,v\frac{\partial\tilde{p}}{\partial\nu}\,dS-\int_D\,v\Delta\tilde{p}\,dx
\\
\\
\displaystyle
&
\displaystyle
=-\int_D\,v(\Delta\tilde{p}+V_0(x)\tilde{p})\,dx+\int_{\partial D}\,v\frac{\partial\tilde{p}}{\partial\nu}\,dS.
\end{array}
$$
Therefore (2.13) becomes
$$\begin{array}{ll}
\displaystyle
\Vert w\Vert_{L^2(\Omega\setminus\overline{D})}^2
&
\displaystyle
=-\int_{\partial D}v\left(\frac{\partial\tilde{p}}{\partial\nu}
+\lambda(x)\tilde{p}\,\right)\,dS
+\int_D\,v(\Delta\tilde{p}+V_0(x)\tilde{p})\,dx
\\
\\
\displaystyle
&
\displaystyle
=\int_D\,v(\Delta\tilde{p}+V_0(x)\tilde{p})\,dx.
\end{array}
$$
This together with (2.12) yields the desired estimate.

\noindent
$\Box$

\subsection{Proof Theorem 1.1 completed}

From Lemma 2.1, (1.19) and the trace theorem one gets the upper bound on (1.20).
The problem is to give the lower bound on (1.20).
We have
$$\begin{array}{l}
\displaystyle
\,\,\,\,\,\,
\left\vert-2\int_{\partial D}\,\text{Im}\,(\lambda(x))\,\text{Im}\,(w\overline{v})\,dS\right\vert
\\
\\
\displaystyle
\le
C\int_{\partial D}\,\vert w\vert\vert v\vert\,dS
\\
\\
\displaystyle
\le C\Vert w\vert_{\partial D}\Vert_{L^2(\partial D)}\Vert v\vert_{\partial D}\Vert_{L^2(\partial D)}
\\
\\
\displaystyle
\le C(\Vert w\vert_{\partial D}\Vert_{L^2(\partial D)}^2+\Vert v\vert_{\partial D}\Vert_{L^2(\partial D)}^2\,)
\\
\\
\displaystyle
\le C(\epsilon_1\Vert\nabla w\Vert_{L^2(\Omega\setminus\overline{D})}^2+\epsilon_1^{-1}\Vert w\Vert_{L^2(\Omega\setminus\overline{D})}^2
+\epsilon_2\Vert\nabla v\Vert_{L^2(D)}^2+\epsilon_2^{-1}\Vert v\Vert_{L^2(D)}^2\,),
\end{array}
$$
where $\epsilon_1$ and $\epsilon_2$ are arbitrary positive numbers with $0<\epsilon_j<1$, $j=1,2$.
And similarly
$$\begin{array}{l}
\displaystyle
\,\,\,\,\,\,
\left\vert
-\int_{\partial D}\,\text{Re}\,(\lambda(x))\,\vert w\vert^2\,dS
+
\int_{\partial D}\,\text{Re}\,(\lambda(x))\,\vert v\vert^2\,dS
\right\vert
\\
\\
\displaystyle
\le C(\epsilon_1\Vert\nabla w\Vert_{L^2(\Omega\setminus\overline{D})}^2+\epsilon_1^{-1}\Vert w\Vert_{L^2(\Omega\setminus\overline{D})}^2
+\epsilon_2\Vert\nabla v\Vert_{L^2(D)}^2+\epsilon_2^{-1}\Vert v\Vert_{L^2(D)}^2\,).
\end{array}
$$
Thus from (1.19) one gets
$$\begin{array}{l}
\displaystyle
\,\,\,\,\,\,
\text{Re}\,\left<\frac{\partial v}{\partial\nu}\vert_{\partial\Omega}-\frac{\partial u}{\partial\nu}\vert_{\partial\Omega},\,\overline{v}\vert_{\partial\Omega}
\right>
\\
\\
\displaystyle
\ge (1-2C\epsilon_1)\Vert \nabla w\Vert_{L^2(\Omega\setminus\overline{D})}^2-C(\epsilon_1)\Vert w\Vert_{L^2(\Omega\setminus\overline{D})}^2
+ (1-2C\epsilon_2)\Vert \nabla v\Vert_{L^2(D)}^2-C(\epsilon_2)\Vert v\Vert_{L^2(D)}^2.
\end{array}
$$
Choosing $\epsilon_1$ and $\epsilon_2$ so small, one obtain
$$
\displaystyle
\text{Re}\,\left<\frac{\partial v}{\partial\nu}\vert_{\partial\Omega}-\frac{\partial u}{\partial\nu}\vert_{\partial\Omega},\,\overline{v}\vert_{\partial\Omega}
\right>
\ge 
C\Vert \nabla v\Vert_{L^2(D)}^2-C'\Vert v\Vert_{L^2(D)}^2-C''\Vert w\Vert_{L^2(\Omega\setminus\overline{D})}^2.
$$
Now applying Lemma 2.2 to this right-hand side, we obtain the lower bound on (1.20).

\noindent
$\Box$

\section{Proof of Theorems 1.2 and 1.3}

\subsection{Two lemmas}

First we give an $L^2$-estimate of the so-called reflected solution in terms of the $L^1$-norm of $v$ 
which is an arbitrary solution of equation (1.7) \footnote{See also Remark 3.1.}.
The proof is essentially the same as that of Lemma 3.1 in \cite{ISurvey} (taken from unpublished manuscript \cite{IEx})
in which the case $V_0(x)\equiv k^2$ is treated.
Here for reader's convenience we present its proof.

\proclaim{\noindent Lemma 3.1.}
Let $v\in H^1(\Omega)$ be an arbitrary complex valued solution of equation (1.7)
and, under Assumption 3 let $u\in H^1(\Omega)$ be the weak solution of (1.21) with $u=v$ on $\partial\Omega$.
Then, there exists a positive constant $C$ independent of $v$ such that 
$$\displaystyle
\Vert u-v\Vert_{L^2(\Omega)}\le C\Vert Vv\Vert_{L^1(\Omega)}.
\tag {3.1}
$$

\endproclaim

{\it\noindent Proof.}
Set $w=u-v$.  By Assumption 3,
one can find the unique weak solution $p\in H^1(\Omega)$ of (1.22) with $F=\overline{w}$.
By (1.23) and $V_0+V_1\in L^{\infty}(\Omega)$, we have $\Delta p\in L^2(\Omega)$ and $\Vert\Delta p\Vert_{L^2(\Omega)}\le C_0\Vert w\Vert_{L^2(\Omega)}$.
The elliptic regularity up to boundary, for the Laplace operator with Dirichlet boundary condition yields $p\in H^2(\Omega)$ and $\Vert p\Vert_{H^2(\Omega)}\le C_1\,\Vert w\Vert_{L^2(\Omega)}$\footnote{
Or use the general result in \cite{GT} Theorem 8.12 on page 186 and (1.23).}.
By the Sobolev imbedding theorem, we have $\Vert p\Vert_{L^{\infty}(\Omega)}\le C_2\Vert p\Vert_{H^2(\Omega)}$.
Thus one gets
$$\displaystyle
\Vert p\Vert_{L^{\infty}(\Omega)}\le C\Vert w\Vert_{L^2(\Omega)}.
\tag {3.2}
$$

It follows from (1.21) that the function $w$ is the weak solution of 
$$\left\{
\begin{array}{ll}
\displaystyle
\Delta w+V_0(x)\,w+V(x)\,w=-V(x)\,v, & x\in\Omega,
\\
\\
\displaystyle
w=0, & x\in\partial\Omega.
\end{array}
\right.
\tag {3.3}
$$
Using this and the weak form of (1.22) with $F=\overline{w}$, we have the expression
$$\begin{array}{ll}
\displaystyle
\Vert w\Vert_{L^2(\Omega)}^2
&
\displaystyle
=-\int_{\Omega}\nabla p\cdot\nabla w\,dx+\int_{\Omega}\,p(V_0(x)+V(x))w\,dx
\\
\\
\displaystyle
&
\displaystyle
=-\int_{\Omega}\,p V(x)\,v\,dx.
\end{array}
$$
This yields
$$\displaystyle
\Vert w\Vert_{L^2(\Omega)}^2
\le
\Vert p\Vert_{L^{\infty}(\Omega)}\,\Vert V\,v\Vert_{L^1(\Omega)}
$$
and thus from (3.2) we obtain the desired estimate.

\noindent
$\Box$

{\bf\noindent Remark 3.1.} By applying the same argument for the derivation of (3.2) to (3.3), we obtain
$$\displaystyle
\Vert u-v\Vert_{L^{\infty}(D)}\le C\Vert Vv\Vert_{L^2(\Omega)}.
\tag {3.4}
$$
If we assume that $D$ is an open set of $\Bbb R^3$ with $\overline{D}\subset\Omega$ instead of $D\subset\Omega$ in Assumption 2, 
one can make use of the interior regularity estimate \cite{GT} Theorem 8.8 on page 183 directly
for the Poisson equation $\Delta w=\tilde{f}\equiv-V_0(x)w-V(x)w-V(x)v$ in $\Omega$
and obtain $\Vert w\Vert_{H^2(D)}\le C\,(\Vert w\Vert_{L^2(\Omega)}+\Vert Vv\Vert_{L^2(\Omega)}\,)$.  By (1.23) in Assumption 3 with $F=-Vv$, this yields $\Vert w\Vert_{H^2(D)}\le C\Vert Vv\Vert_{L^2(\Omega)}$ and the Sobolev imbedding yields (3.4).
This derivation has an advantage since only the interior regularity result without specifying the boundary condition on $\partial\Omega$ is required at the cost of the strong assumption $D$ being open and satisfying $\overline D\subset\Omega$.  See also Remark 3.2 for an explanation about the role of (3.4) in the proof of Theorems 1.2 and 1.3.

Next we describe another lemma which is originally stated in the unpublished manuscript \cite{IEx} and Lemma 3.2 in \cite{ISurvey}
provided $D$ is an open set of $\Bbb R^3$ and $\partial D\in C^2$ and $\omega$ satisfies:
$$\begin{array}{ll}
\displaystyle
x\cdot\omega<h_D(\omega) & \text{a.e.} x\in \partial D.
\end{array}
$$
The proof presented here is completely different from the original one and removes such assumption.

\proclaim{\noindent Lemma 3.2.}  Assume that $D$ is $p$-regular with respect to $\omega$.
We have
$$\displaystyle
\lim_{\tau\rightarrow\infty}\,\frac{\Vert v_0\Vert_{L^1(D)}}{\Vert v_0\Vert_{L^2(D)}}=0,
\tag {3.5}
$$
where $v_0(x)=e^{\tau x\cdot\omega}$.
\endproclaim
{\it\noindent Proof.} It suffices to prove the formula
$$\displaystyle
\lim_{\tau\rightarrow\infty}\,
\frac{e^{-\tau h_D(\omega)}\,\Vert v_0\Vert_{L^1(D)}}
{e^{-\tau h_D(\omega)}\,\Vert v_0\Vert_{L^2(D)}}
=0.
$$
Let $\epsilon$ be an arbitrary positive number.
We have
$$\begin{array}{ll}
\displaystyle
e^{-\tau h_D(\omega)}\,\Vert v_0\Vert_{L^1(D)}
&
\displaystyle
=\int_D e^{-\tau(h_D(\omega)-x\cdot\omega)}\,dx
\\
\\
\displaystyle
&
\displaystyle
=\int_{D_{\omega}(\epsilon)} e^{-\tau(h_D(\omega)-x\cdot\omega)}\,dx
\\
\\
\displaystyle
&
\displaystyle
\,\,\,
+
\int_{D\setminus D_{\omega}(\epsilon)}e^{-\tau(h_D(\omega)-x\cdot\omega)}\,dx
\\
\\
\displaystyle
&
\displaystyle
\le
\left(\int_{D_{\omega}(\epsilon)}dx\right)^{\frac{1}{2}}\left(\int_De^{-2\tau(h_D(\omega)-x\cdot\omega)}\,dx\right)^{\frac{1}{2}}
\\
\\
\displaystyle
&
\displaystyle
\,\,\,
+\int_D\,dx e^{-\tau\epsilon}.
\end{array}
\tag {3.6}
$$
Since $D_{\omega}(\epsilon)$ is contained in the set $\Omega\cap\{x\in\Bbb R^3\,\vert\,h_D(\omega)-\epsilon<x\cdot\omega\le h_D(\omega)\}$,
we have
$$\displaystyle
\left(\int_{D_{\omega}(\epsilon)}dx\right)^{\frac{1}{2}}
\le C_1\epsilon^{\frac{1}{2}}.
$$
Thus (3.6) yields
$$
\displaystyle
e^{-\tau h_D(\omega)}\,\Vert v_0\Vert_{L^1(D)}
\le C_2\left(\epsilon^{\frac{1}{2}}e^{-\tau h_D(\omega)}\Vert v\Vert_{L^2(D)}
+e^{-\tau\epsilon}
\right).
\tag {3.7}
$$ 
By the $p$-regularity of $D$ with respect to $\omega$ and Fubini's theorem, for all $\tau>>1$ we have
$$\begin{array}{ll}
\displaystyle
\int_{D}e^{-2\tau (h_D(\omega)-x\cdot\omega)}\,dx
&
\displaystyle
\ge \int_{0}^{\delta}\mu_2(S_{\omega}(s))e^{-2\tau s}\,ds
\\
\\
\displaystyle
&
\displaystyle
\ge C\int_{0}^{\delta}s^{p-1}e^{-2\tau s}\,ds
\\
\\
\displaystyle
&
\displaystyle
=\frac{C}{\tau^{p}}\,\int_0^{\tau\delta}s^{p-1}e^{-2s}\,ds
\\
\\
\displaystyle
&
\displaystyle
\ge \frac{C_3}{\tau^{p}}.
\end{array}
$$
Thus one gets
$$\displaystyle
e^{-\tau h_D(\omega)}\Vert v_0\Vert_{L^2(D)}\ge C_4\tau^{-\frac{p}{2}}.
\tag {3.8}
$$
This together with (3.7) yields
$$\displaystyle
\frac{e^{-\tau h_D(\omega)}\,\Vert v_0\Vert_{L^1(D)}}
{e^{-\tau h_D(\omega)}\,\Vert v_0\Vert_{L^2(D)}}
\le
C_4
\left(\epsilon^{\frac{1}{2}}+e^{-\tau\epsilon}\tau^{\frac{p}{2}}\,\right).
\tag {3.9}
$$
Now set $\epsilon=\tau^{-\gamma}$ with an arbitrary fixed number $\gamma\in\,]0,\,1[$.
Then as $\tau\rightarrow\infty$ we have
$$\begin{array}{ll}
\displaystyle
\epsilon^{\frac{1}{2}}+e^{-\tau\epsilon}\tau^{\frac{p}{2}}
&
\displaystyle
=\tau^{-\frac{\gamma}{2}}+e^{-\tau^{1-\gamma}}\tau^{\frac{p}{2}}
\\
\\
\displaystyle
&
\displaystyle
\rightarrow 0.
\end{array}
$$
A combination of this and (3.9) yields (3.5).

\noindent
$\Box$

If $D$ is open and $\partial D$ is Lipschitz, then $D$ is $3$-regular with respect to all unit vectors $\omega$.  Thus
Lemma 3.2 is valid also for such $D$ without any restriction on $\omega$.  The advantage of the proof above is: we do not make use
of any local coordinate system of $D$ in a neighbourhood of $\partial D$.

\subsection{Finishing the proof}

Let us continue to prove Theorem 1.2.
We give a proof only the case when $\text{Re}\,V$ has a positive jump since another case can be treated similarly.

By the Alessandrini identity and Assumption 3, from (1.32) we have
$$\begin{array}{ll}
\displaystyle
I_{\omega,\vartheta}(\tau)
&
\displaystyle
=\int_{D}\,V(x)vv^*\,dx+\int_{D}\,V(x)\,wv^*\,dx,
\end{array}
\tag {3.10}
$$
where $w=u-v$.

Write
$$\begin{array}{ll}
\displaystyle
\int_{D}\,V(x)vv^*\,dx
&
\displaystyle
=\int_{D}\,V(x)e^{2\tau x\cdot\omega}(1+\Psi(x,z))(1+\Psi(x,\overline{z}))\,dx
\\
\\
\displaystyle
&
\displaystyle
=\int_{D}\,V(x)e^{2\tau x\cdot\omega}\,dx
\\
\\
\displaystyle
&
\displaystyle
\,\,\,
+\int_{D}\,V(x)e^{2\tau x\cdot\omega}(\Psi(x,z)+\Psi(x,\overline{z})+\Psi(x,z)\,\Psi(x,\overline{z}))\,dx.
\end{array}
$$
By (3.1) and the Ramm estimate $\Vert\Psi(\,\cdot\,,z)\Vert_{L^{\infty}(\Omega)}=\Vert\Psi(\,\cdot\,,\overline{z})\Vert_{L^{\infty}(\Omega)}=\epsilon(\tau^{-1})$,
we have
$$\displaystyle
\left\vert
\int_{D}\,V(x)e^{2\tau x\cdot\omega}(\Psi(x,z)+\Psi(x,\overline{z})+\Psi(x,z)\,\Psi(x,\overline{z}))\,dx
\right\vert
=\Vert v_0\Vert_{L^2(D)}^2\epsilon(\tau^{-1}),
$$
where $v_0(x)=e^{\tau x\cdot\omega}$.
This together with (1.27), (1.29), (3.1) and (3.10) gives
$$
\displaystyle
\text{Re}\,I_{\omega,\vartheta}(\tau)
\ge \int_D\,\text{Re}\,V(x)\vert v_0\vert^2\,dx-\Vert v_0\Vert_{L^2(D)}^2\epsilon(\tau^{-1})
-C\Vert v_0\Vert_{L^2(D)}\Vert v_0\Vert_{L^1(D)}.
$$

Rewrite this as
$$\begin{array}{l}
\displaystyle
\,\,\,\,\,\,
e^{-2\tau h_D(\omega)}\text{Re}\,I_{\omega,\vartheta}(\tau)
\\
\\
\displaystyle
\ge
e^{-2\tau h_D(\omega)}
\,
\Vert v_0\Vert_{L^2(D)}^2
\,\left\{
\frac{\displaystyle
\int_D\,\text{Re}\,V(x)\vert e^{-\tau h_D(\omega)}\,v_0\vert^2\,dx}
{\displaystyle
e^{-2\tau h_D(\omega)}\,\Vert v_0\Vert_{L^2(D)}^2}
-\epsilon(\tau^{-1})
-C\,\frac{\Vert v_0\Vert_{L^1(D)}}{\Vert v_0\Vert_{L^2(D)}}\,
\right\}.
\end{array}
\tag {3.11}
$$

Since we have $\text{Re}\,V(x)\ge C$ for almost all $x\in D_{\omega}(\delta)$,
we see that
$$\begin{array}{l}
\displaystyle
\,\,\,\,\,\,
\int_D\,\text{Re}\,V(x)\,\vert e^{-\tau h_D(\omega)}\,v_0\vert^2\,dx
\\
\\
\displaystyle
\displaystyle
=\int_{D_{\omega}(\delta)}\,\text{Re}\,V(x)\vert e^{-\tau h_D(\omega)}\,v_0\vert^2\,dx
+\int_{D\setminus D_{\omega}(\delta)}\,\text{Re}\,V(x)\vert e^{-\tau h_D(\omega)}\,v_0\vert^2\,dx
\\
\\
\displaystyle
\ge
C\int_{D_{\omega}(\delta)}\,\vert e^{-\tau h_D(\omega)}\,v_0\vert^2\,dx
-C'e^{-2\tau\delta}
\\
\\
\displaystyle
=Ce^{-2\tau h_D(\omega)}\,\Vert v_0\Vert_{L^2(D)}^2
-
C\int_{D\setminus D_{\omega}(\delta)}\,\vert e^{-\tau h_D(\omega)}\,v_0\vert^2\,dx
+O(e^{-2\tau\delta})
\\
\\
\displaystyle
=Ce^{-2\tau h_D(\omega)}\,\Vert v_0\Vert_{L^2(D)}^2
+O(e^{-2\tau\delta}).
\end{array}
$$
Thus this together with (3.8) and (3.11) yields:
there exist $C_1>0$ and $\tau_1>0$ such that, for all $\tau\ge\tau_1$
$$\displaystyle
\frac{\displaystyle
\int_D\,\text{Re}\,V(x)\,\vert e^{-\tau h_D(\omega)}\,v\vert^2\,dx}
{e^{-2\tau h_D(\omega)}\,\Vert v\Vert_{L^2(D)}^2}
\ge C_1.
$$
Now from this together with (3.8) again, (3.11) and (3.5) we conclude that 
there exist $C_2>0$ and $\tau_2>0$ such that, for all $\tau\ge\tau_2$ 
$$\displaystyle
e^{-2\tau h_D(\omega)}\tau^p\text{Re}\,I_{\omega,\vartheta}(\tau)\ge C_2.
$$

On the other hand, from Lemma 3.1 combined with $\Vert Vv\Vert_{L^1(D)}\le C\Vert v\Vert_{L^2(D)}$ and the Sylester-Uhlmann estimate
restricted to $L^2(\Omega)$,  we see that $e^{-2\tau h_D(\omega)}\,\vert I_{\omega,\vartheta}(\tau)\vert$ is bounded 
from above for all $\tau>>1$.  Now a standard argument yields the desired conclusions.

The proof of Theorem 1.3 can be done in the same way, so its description is omitted.

$\quad$

{\bf\noindent Remark 3.2.}
The point of the derivation of (3.11) is the estimate of the following type
$$\displaystyle
\left\vert\int_{\Omega}Q(x)(u-v)v^*\,dx\right\vert\le C\Vert v_0\Vert_{L^2(D)}\,\Vert v_0\Vert_{L^1(D)},
$$
where $Q\in L^{\infty}(\Omega)$ and satisfies $Q(x)=0$ a. e. $x\in\Omega\setminus D$.  
If $\overline{D}\subset\Omega$ instead of $D\subset\Omega$, one can use the local interior estimate (3.4).
In fact, we have
$$\begin{array}{ll}
\displaystyle
\left\vert\int_{\Omega}Q(x)(u-v)v^*\,dx\right\vert
&
\displaystyle
\le\Vert Q\Vert_{L^{\infty}(D)}\Vert u-v\Vert_{L^{\infty}(D)}\Vert v^*\Vert_{L^1(D)}
\\
\\
\displaystyle
&
\displaystyle
\le C\Vert Q\Vert_{L^{\infty}(D)}\Vert Vv\Vert_{L^2(D)}\Vert v_0\Vert_{L^1(D)}
\\
\\
\displaystyle
&
\displaystyle
\le C'\Vert v_0\Vert_{L^2(D)}\Vert v_0\Vert_{L^1(D)}.
\end{array}
$$
Note that the upper estimate on (1.27) is also used.

\section{Proof of Lemma 3.2 via Sini-Yoshida's approach}

First we present a {\it formal} application of their approach in the proof of a lemma in \cite{SY} to the proof of Lemma 3.2
in the case when $\partial D$ is Lipschitz.
Second the author would like to point a problem out whether their idea really can cover the case when $\partial D$ is {\it Lipschitz}.

Given a sufficiently small positive number $\delta$ they cover the compact set $\{x\in\partial D\,\vert\,x\cdot\omega=h_D(\omega)\}$ 
by a finite set $\{C_j\}_{j=1,\cdots,N}$ of the rotation of around point $x_j\in\partial D$, $j=1,\cdots,N$ with $x_j\cdot\omega=h_D(\omega)$ and translation to $x_j$
of a cubic domain centered at the origin $(0,0,0)$ which corresponds to point $x_j$
\footnote{In fact, instead of cubic, they use open ball, however, there is no essential difference.} in $\Bbb R^3$
in such a way that

$\bullet$ each $C_j\cap D$ has the expression 
$$\displaystyle
C_j\cap D=\{x=x_j+A_jy\,\vert\,l_j(y')<y_3<\delta,\,\vert y'\vert<\delta\,\},
\tag {4.1}
$$
where $y'=(y_1,y_2)$, $l_j(y')$ is a Lipschitz continuous function on $\Bbb R^2$  and satisfies $0\le l_j(y')\le C\vert y'\vert$
with a positive constant $C$ and $A_j$ is an orthogonal matrix.

If necessary, choosing a smaller $\delta$ one may assume that
$$\displaystyle
C_{\delta}\equiv\inf_{x\in D\setminus(\cup_{j=1}^N\,C_j\cap D)}\text{dist}\,(\{x\},\pi_{\omega})>0,
$$
where $\pi_{\omega}$
denotes the plane $x\cdot\omega=h_D(\omega)$.

The point is the choice of $A_j$ in each $C_j\cap D$.
For $x=x_j+A_jy$, we have
$$\begin{array}{ll}
\displaystyle
h_D(\omega)-x\cdot\omega
&
\displaystyle
=(x_j-x)\cdot\omega
\\
\\
\displaystyle
&
\displaystyle
=-A_jy\cdot\omega
\\
\\
\displaystyle
&
\displaystyle
=-y\cdot A_j^T\omega.
\end{array}
$$
From their computation in\cite{SY}, it is clear that they choose each $A_j$ in such a way that
$$\displaystyle
A_j^T\omega=-\mbox{\boldmath $e_3$}.
\tag {4.2}
$$
Thus we have
$$\begin{array}{ll}
\displaystyle
e^{-\tau h_D(\omega)}v_0
=e^{-\tau y_3}, & x=x_j+A_jy\in C_j\cap D.
\end{array}
\tag {4.3}
$$

Since $C_j\cap D\subset D$, $j=1,\cdots, N$ one has the {\it lower and upper bound}:
$$\displaystyle
\frac{1}{N}M(\tau)
\le
e^{-\tau h_D(\omega)}\int_D e^{\tau x\cdot\omega}\,dx
\le 
M(\tau)
+O(e^{-C_{\delta}\tau}),
\tag {4.4}
$$
where 
$$\displaystyle
M(\tau)
=\sum_{j=1}^Ne^{-\tau h_D(\omega)}\int_{C_j\cap D}\,e^{\tau x\cdot\omega}\,dx,
$$
The point is, in both the lower and upper bound the same estimator $M(\tau)$
satisfying $M(2\tau)\le C M(\tau)$ appears.   Besides by using the local expression of $C_j\cap D$ and (4.3) the each term in $M(\tau)$ 
one has
$$\begin{array}{ll}
\displaystyle
M(\tau) 
&
\displaystyle
=\sum_{j=1}^N\int_{\vert y'\vert<\delta}\,dy'\int_{l_j(y')}^{\delta}e^{-\tau y_3}\,dy_3
\\
\\
\displaystyle
&
\displaystyle
=\tau^{-1}\sum_{j=1}^N\int_{\vert y'\vert<\delta} e^{-\tau l_j(y')}\,dy'+O(e^{-C_{\delta}'\tau}),
\end{array}
$$
where $l_j(y')$ is a Lipschitz continuous function and satisfying $0\le 0\le l_j(y')\le C\vert y'\vert$
with a positive constant $C$.  Note that this yields a rough lower estimate
$$\displaystyle
M(2\tau)\ge C\tau^{-3}.
\tag {4.5}
$$

By virtue of the Schwarz inequality the $M(\tau)$ has the property:
$$\begin{array}{ll}
\displaystyle
(M(\tau))^2
&
\displaystyle
=O(\tau^{-2}\left(\sum_{j=1}^N\int\int_{\vert y'\vert<\delta} e^{-\tau l_j(y')}\,dy'\right)^2+e^{-2C_{\delta'}\tau})
\\
\\
\displaystyle
&
\displaystyle
=O(\tau^{-2}\sum_{j=1}^N\int\int_{\vert y'\vert<\delta}e^{-2\tau l_j(y')}\,dy'+e^{-2C_{\delta}'\tau})
\\
\\
\displaystyle
&
\displaystyle
=O(\tau^{-1}M(2\tau)+e^{-2C_{\delta}''\tau}).
\end{array}
$$
A combination of this, the right-hand side on (4.4) and (4.5) gives
$$\begin{array}{ll}
\displaystyle
\frac{e^{-2\tau h_D(\omega)}\,\Vert v_0\Vert_{L^1(D)}^2}
{e^{-2\tau h_D(\omega)}\,\Vert v_0\Vert_{L^2(D)}^2}
&
\displaystyle
\le
\frac{(M(\tau)+O(e^{-\tau C_{\delta}}))^2}
{\frac{1}{N}M(2\tau)}\\
\\
\displaystyle
&
\displaystyle
\le
C\frac{M(\tau)^2+O(e^{-2\tau C_{\delta}})}
{M(2\tau)}
\\
\\
\displaystyle
&
\displaystyle
\le C\frac{\tau^{-1}M(2\tau)+e^{-2C_{\delta}'''\tau}}{M(2\tau)}
\\
\\
\displaystyle
&
\displaystyle
=C\left(\tau^{-1}+\frac{e^{-2C_{\delta}'''\tau}}{M(2\tau)}\right)
\\
\\
\displaystyle
&
\displaystyle
=O(\tau^{-1}).
\end{array}
$$
Thus one gets (3.5) with convergence rate is $O(\tau^{-\frac{1}{2}})$.

The author thinks that their idea of proof has an advantage since it never makes use of any {\it concrete} upper bound of $M(\tau)$ as $\tau\rightarrow\infty$
which needs some additional condition about the behaviour of $l_j(y')$ as $y'\rightarrow (0,0)$.

\subsection{What is a problem on their approach?}

The problem is the choice of $A_j$ in such a way that $A_j$ satisfies (4.2) and  $C_j\cap D$ has the expression (4.1) at the same time.

Consider, for simplicity, two dimensional case.
Let $D$ be the domain defined by
$$
D=\{x\in\Bbb R^2\,\vert\,x_1<x_2<2x_1\,,0<x_1<\infty\,\}\cap B,
$$
where $B$ is a unit open disk centered the origin $(0,0)$.  The $D$ is Lipschitz \cite{Gr}.
Let $\omega=-\mbox{\boldmath $e$}_2$.
Then we have 
$$\displaystyle
\partial D\cap\pi_{\omega}=\{(0,0)\}.
$$
Let $C$ be the square given by
$$\displaystyle
C=\{y\in\Bbb R^2\, \vert\, \vert y_1\vert<\delta,\,\vert y_2\vert<\delta\,\}
$$
with an arbitrary fixed positive number $\delta\in\,]0,\,\frac{1}{\sqrt{2}}\,[$.  Let $A(\theta)$ be a rotation around $(0,0)$ given by the matrix
$$\begin{array}{ll}
\displaystyle
A(\theta)=\left(\begin{array}{ll}
\displaystyle
\cos\theta & \displaystyle -\sin\theta\\
\\
\displaystyle
\sin\theta & \displaystyle\cos\theta
\end{array}
\right), & \theta\in\,[0,\,2\pi[.
\end{array}
$$
Set $C(\theta)=A(\theta)C$ which is the rotation of $C$ around $(0,0)\in\partial D\cap\pi_{\omega}$.
Then we have
$$\displaystyle
C(\theta)\cap D=\{x=A(\theta)y\,\vert\,y\in C\cap\hat{D}(\theta)\,\},
$$
where
$$\displaystyle
\hat{D}(\theta)
=\{(y_1,y_2)\in C\,\vert\,y_1(\cos\theta-\sin\theta)<y_2(\cos\theta+\sin\theta),\,
y_2(\cos\theta+2\sin\theta)<y_1(2\cos\theta-\sin\theta)\,\}.
$$

Since $\omega=-\mbox{\boldmath $e$}_2$ and $h_D(\omega)=0$, by the change of variables $x=A(\theta)y$ we have
$$\begin{array}{ll}
\displaystyle
h_D(\omega)-x\cdot\omega
&
\displaystyle
=-A(\theta)y\cdot\omega
\\
\\
\displaystyle
&
\displaystyle
=y\cdot A(-\theta)\mbox{\boldmath $e$}_2\\
\\
\displaystyle
&
\displaystyle
=y_1\sin\theta+y_2\cos\theta.
\end{array}
$$
This yields the expression
$$\begin{array}{ll}
\displaystyle
e^{-\tau h_D(\omega)}v_0
=e^{-\tau (y_1\sin\theta+y_2\cos\theta)}, & x=A(\theta)y\in C(\theta)\cap D.
\end{array}
\tag {4.6}
$$
Thus this expression becomes (4.3) type, that is $e^{-\tau(h_D(\omega)-x\cdot\omega)}=e^{-\tau y_2}$ if and only if
$\theta=0$.  Then we see that $C(0)\cap D=D$ never coincides with a domain having the form of the two dimensional version of (4.1).
Of course, by choosing another $\theta$ one can make $C(\theta)\cap D$ such as the domain of (4.1) type, however, in that case
one can not obtain the expression $e^{-\tau(h_D(\omega)-x\cdot\omega)}=e^{-\tau y_2}$ for $x=A(\theta)\in C(\theta)\cap D$.

In three-dimensional case, consider a domain having a cone or convex polyhedron as a part, then one encounters a more complicated situation.
Therefore in the Lipschitz case, their approach does not work as they wished.
In this sense our approach is better to cover a broad class of boundary.  
Note that the comment also works for the articles \cite{KLS} and \cite{KS} following their paper.

\section{Applications}

\subsection{Revisiting the enclosure method}

First we consider the simplest case.  Let $k$ be a fixed positive number.
We identify the set of all real unit vectors with the unit sphere $S^2$.
Given an arbitrary $\omega\in S^2$ and $\vartheta\in S^2$ with $\omega\cdot\vartheta=0$ define
$$\begin{array}{ll}
\displaystyle
v(x,z)=e^{x\cdot z}, & x\in\Bbb R^3
\end{array}
$$
where $z=\tau\omega+i\sqrt{\tau^2+k^2}\vartheta$ and $\tau>0$.  Since the complex vector $z$ satisfies the equation $z\cdot z=-k^2$,
the $v=v(\,\,\cdot\,;z)$ is an entire solution of the Helmholtz equation $\Delta v+k^2v=0$.

Define the indicator function of the enclosure method by the formula:
$$\begin{array}{ll}
\displaystyle
I_{\omega,\vartheta}(\tau)=
\text{Re}\,\left<\frac{\partial v}{\partial\nu}\vert_{\partial\Omega}
-\frac{\partial u}{\partial\nu}\vert_{\partial\Omega},\,\overline{v}\vert_{\partial\Omega}
\right>, & \tau>0,
\end{array}
$$
where $v=v(\,\cdot\,;z)$ and $u$ is the weak solution of (1.10) with $V_0(x)\equiv k^2$.
The following result is an extension of the author's previous result Theorem E in \cite{Iimpedance}.

\proclaim{\noindent Theorem 5.1.}  Assume that $\partial D$ is $C^{1,1}$ and $\lambda\in C^{0,1}(\partial D)$.  
Under Assumptions 1 with $V_0(x)\equiv k^2$, 
we have, for all $\gamma\in\Bbb R$
$$\displaystyle
\lim_{\tau\rightarrow\infty}\,e^{-2\tau t}\tau^{\gamma}\,I_{\omega,\vartheta}(\tau)
=
\left\{
\begin{array}{ll}
\displaystyle
0 & \text{if $t>h_D(\omega)$,}
\\
\\
\displaystyle
\infty & \text{if $t<h_D(\omega)$,}
\end{array}
\right.
$$
and
$$\begin{array}{ll}
\displaystyle
\liminf_{\tau\rightarrow\infty}e^{-2\tau t}\tau\,I_{\omega,\vartheta}(\tau)>0 & \text{if $t=h_D(\omega)$.}
\end{array}
$$
Besides, we have the one-line formula
$$\displaystyle
\lim_{\tau\rightarrow\infty}\,\frac{1}{2\tau}\,\log\,I_{\omega,\vartheta}(\tau)=h_D(\omega).
$$
\endproclaim
{\it\noindent Proof.}
Theorem 1.1 gives
$$\displaystyle
C_1\Vert\nabla v\Vert_{L^2(D)}^2-C_2\Vert v\Vert_{L^2(D)}^2
\le 
I_{\omega,\vartheta}(\tau)
\le C_3\Vert v\Vert_{H^1(D)}^2.
\tag {5.1}
$$
Since we have $\Vert\nabla v\Vert_{L^2(D)}^2=(2\tau^2+k^2)\Vert v\Vert_{L^2(D)}^2$, (5.1) yields, for all $\tau>>1$
$$\displaystyle
C_4\tau^2e^{-2\tau h_D(\omega)}\Vert v\Vert_{L^2(D)}^2\le e^{-2\tau h_D(\omega)}I_{\omega,\vartheta}(\tau)\le C_5\tau^2e^{-2\tau h_D(\omega)}\Vert v\Vert_{L^2(D)}^2.
$$
Since $\partial D$ is $C^{1,1}$ and thus Lipschitz.  Then $\partial D$ is $3$-regular with respect to all $\omega\in S^2$
and this yields\footnote{Of course this is a rough lower estimate.}$e^{-2\tau h_D(\omega)}\Vert v\Vert_{L^2(D)}^2\ge C_6\tau^{-3}$.
Noting the trivial estimate $e^{-2\tau h_D(\omega)}\Vert v\Vert_{L^2(D)}\le C_7$, we finally obtain
$$\displaystyle
C_8e^{2\tau h_D(\omega)}\le \tau I_{\omega,\vartheta}(\tau)\le C_9\tau^3 e^{2\tau h_D(\omega)}.
$$ 
Now from this all the assertions of Theorem 5.1 are valid.

\noindent
$\Box$

In \cite{SY} they proved a corresponding result to Theorem 5.1 for the case when $\lambda=0$ and $\partial D$ is Lipschitz.
Needless to say it is almost a direct consequence of (1.3) as we have seen above
\footnote{Strangely enough after having (1.3) their proof does not make use of 
the explicit relation between $v$ and $\nabla v$.  See line 14 down on page 15 in \cite{SY} and the proof of Lemma 3.7 therein.}.

Next consider the general case.  Let $V_0\in C^{0,1}(\overline{\Omega})$.
A combination of Theorem 6.2.4 on p.277 in \cite{Gr} 
and a cut-off function, one has a $\tilde{V}_0\in C^{0,1}(\Bbb R^3)$ with compact support such that
$V_0(x)=\tilde{V_0}(x)$ for all $x\in\overline{\Omega}$.  
Replace $\tilde{V_0}$ in equation (1.25) by this one and solve it.
Using the solution, we define $v=v(\,\cdot\,,z)$ with $\tau>>1$ be the same form as (1.26).
Since the new $\tilde{V_0}$ belongs to $H^1(\Bbb R^3)$ with compact support, from the Ramm estimate we have 
$$\displaystyle
\lim_{\tau\rightarrow\infty}\Vert\nabla\Psi(\,\cdot\,,z)\vert_{\Omega}\Vert_{L^{\infty}(\Omega)}=0.
$$
This implies that not only estimate (1.27) for $v$ but also for $\nabla v$ for all $\tau>>1$:
$$\begin{array}{ll}
\displaystyle
C_1\tau e^{x\cdot\omega}\le\vert \nabla v(x,z)\vert\le C_2 \tau e^{x\cdot\omega}, & \text{a.e.}x\in\Omega.
\end{array}
\tag {5.2}
$$
Let $u$ be the weak solution of (1.10) with $v=v(\,\cdot\,,z)$ and 
define the indicator function $I_{\omega,\vartheta}(\tau)$ as the exactly same form as above.
Then as a direct corollary of Theorem 1.1 combined with (1.27) and (5.2) we immediately obtain the following result.

\proclaim{\noindent Theorem 5.2.}  Assume that $\partial D$ is $C^{1,1}$, $\lambda\in C^{0,1}(\partial D)$
and $V_0\in C^{0,1}(\overline{\Omega})$.  Under Assumptions 1-2, we have the same conclusions as Theorem 5.1.
\endproclaim

\subsection{A simple proof of the side A of the probe method}

In this subsection it is assumed that $V_0(x)\equiv k^2$ with a fixed positive number $k$ satisfies Assumption 1, $\partial D$ is $C^{1,1}$ and
$\lambda\in C^{0,1}(\partial D)$.
By a standard lifting argument, given $y\in\Omega\setminus\overline{D}$ there exists the unique weak solution
$w_y\in H^1(\Omega\setminus\overline{D})$
of
$$\left\{
\begin{array}{ll}
\displaystyle
\Delta w_y+k^2w_y=0, & x\in\Omega\setminus\overline D,\\
\\
\displaystyle
\frac{\partial w_y}{\partial\nu}+\lambda(x)\,w_y=-\frac{\partial G_y}{\partial\nu}-\lambda(x)\,G_y
, & x\in\partial D,
\\
\\
\displaystyle
w_y=0, & x\in\partial\Omega.
\end{array}
\right.
$$
where
$$\begin{array}{ll}
\displaystyle
G_y(x)=\frac{e^{ik\vert x-y\vert}}{4\pi\vert x-y\vert}, & x\not=y.
\end{array}
$$

Recall the indicator function in the side A of the probe method, Definition 2.3 in \cite{Iimpedance}:
$$\begin{array}{ll}
\displaystyle
\,\,\,\,\,\,
I(y)
&
\displaystyle
=-2\int_{\partial D}\,\text{Im}\,(\lambda(x))\,\text{Im}\,(w_y\overline{G_y})\,dS
\\
\\
\displaystyle
&
\displaystyle
\,\,\,
-\int_{\partial D}\,\text{Re}\,(\lambda(x))\,\vert w_y\vert^2\,dS+\int_{\Omega\setminus\overline{D}}\,\vert\nabla w_y\vert^2\,dx-\int_{\Omega\setminus\overline{D}}\,\text{Re}\,(V_0(x))\,\vert w_y\vert^2\,dx
\\
\\
\displaystyle
&
\displaystyle
\,\,\,
+\int_{\partial D}\,\text{Re}\,(\lambda(x))\,\vert G_y\vert^2\,dS+\int_D\,\vert\nabla G_y\vert^2\,dx-\int_D\,\text{Re}\,(V_0(x))\,\vert G_y\vert^2\,dx.
\end{array}
$$

Given an arbitrary point $y\in\Omega$ a needle $\sigma$ with a tip at $y$ is a non self-intersecting piecewise linear curve with a parameter $t\in\,[0,\,1]$ 
such that $\sigma(0)\in\partial\Omega$, $\sigma(1)=y$, and $\sigma(t)\in\Omega$ for all $t\in\,]0,\,1[$.
Let $\xi=\{v_n\}$ be an arbitrary needle sequence for $(y,\sigma)$,
that is, each $v_n\in H^1(\Omega)$ is a weak solution of the Helmholtz equation in $\Omega$
and, for all compact set $K$ of $\Bbb R^3$ satisfying $K\subset\Omega\setminus\sigma$ it holds that
$$\displaystyle
\lim_{n\rightarrow\infty}(\Vert v_n-G_y\Vert_{L^2(K)}+\Vert\nabla(v_n-G_y)\Vert_{L^2(K)})=0.
$$
It is known at the early stage of the probe method \cite{IProbe2}(see also \cite{IProbeNew}) that under the additional assumption
that $k^2$ is not a Dirichlet eigenvalue of $-\Delta$ in $\Omega$, there exists a needle sequence for an arbitrary $(y,\sigma)$ by the Runge approximation property of the Helmholtz equation.
Besides, if $\sigma$ is given by a part of a line which is called a straight needle, 
then one can give an explicit needle sequence without any restriction on $k^2$, see \cite{ICarleman}.

Recall also the indicator sequence, Definition 2.2 in \cite{Iimpedance}:

$$\begin{array}{ll}
\displaystyle
I(y,\sigma,\xi)_n
&
\displaystyle
=\text{Re}\,\left<\frac{\partial v_n}{\partial\nu}\vert_{\partial\Omega}-\frac{\partial u_n}{\partial\nu}\vert_{\partial\Omega}, \overline{v_n}\,\right>,
\end{array}
$$
where $\xi=\{v_n\}$ is a needle sequence for $(y,\sigma)$ and $u_n$ is the weak solution of 
$$\left\{
\begin{array}{ll}
\displaystyle
\Delta u_n+k^2u_n=0, & x\in\partial\Omega\setminus\overline D,\\
\\
\displaystyle
\frac{\partial u_n}{\partial\nu}+\lambda(x)\,u_n=0, & x\in\partial D,
\\
\\
\displaystyle
u_n=v_n, & x\in\partial\Omega.
\end{array}
\right.
$$
The following result gives us (i) a computation procedure of the indicator function from the indicator sequence; (ii) qualitative behaviour
of the computed indicator function at the places away from and near obstacle.
It is an extension of the result in \cite{IProbe2} when $\partial D$ is $C^2$ and  $\lambda(x)\equiv 0$
and the result in \cite{Iimpedance} when $\partial D\in C^2$ and $\lambda\in C^1(\partial D)$ with 
$\min_{x\in\partial D}\text{Im}\,\lambda(x)>0$\footnote{This is a sufficient condition that ensures the validity of Assumption 1.}
to the case $\partial D\in C^{1,1}$ and $\lambda\in C^{0,1}(\partial D)$.

\proclaim{\noindent Theorem 5.3.}  Let $k$ be an arbitrary fixed positive number satisfying Assumption 1 with $V_0(x)\equiv k^2$.
Assume that $\partial D$ is $C^{1,1}$ and $\lambda\in C^{0,1}(\partial D)$.

(a)  Given $y\in\Omega\setminus\overline D$ we have 
$$\displaystyle
I(y)=\lim_{n\rightarrow\infty}\,I(y,\sigma,\xi)_n,
$$
where $\sigma$ is an arbitrary needle with tip at $y$ such that $\sigma\cap\overline{D}=\emptyset$ and
$\xi=\{v_n\}$ is an arbitrary needle sequence for $(y,\sigma)$.

(b)  For each $\epsilon>0$ we have
$$\displaystyle
\sup_{\text{dist}\,(y,D)>\epsilon}\,I(y)<\infty.
$$

(c)  For any point $a\in\partial D$ we have
$$\displaystyle
\lim_{y\rightarrow a}\,I(y)=\infty.
$$

\endproclaim
{\it\noindent Proof.}  We have $v_n\rightarrow G_y$ in $H^1(D)$ provided $\sigma\cap\overline{D}=\emptyset$.
This together with Lemma 2.1 yields $w_n\equiv u_n-v_n\rightarrow w_y$ in $H^1(\Omega\setminus\overline{D})$.
Thus (1.19) yields the validity of (a).
Then from Theorem 1.1 one gets
$$\displaystyle
C_1\Vert\nabla G_y\Vert_{L^2(D)}^2-C_2\Vert G_y\Vert_{L^2(D)}^2
\le I(y)\le C_3\Vert G_y\Vert_{H^1(D)}^2.
$$
The statements (b) and (c) are direct consequence of this
since we have the estimates
$$\begin{array}{l}
\displaystyle
\sup_{y\in\Omega}\Vert G_y\Vert_{L^2(D)}<\infty,
\\
\\
\displaystyle
\sup_{y\in\Omega,\,\text{dist}\,(y,D)>\epsilon}\Vert G_y\Vert_{H^1(D)}<\infty
\end{array}
$$
and
$$\displaystyle
\lim_{y\rightarrow a,\,y\in\Omega\setminus\overline{D}}\,\Vert\nabla G_y\Vert_{L^2(D)}=\infty.
$$

\noindent
$\Box$

It should be emphasized that, by virtue of Theorem 1.1 the proof given here is so simple compared with the previous one for impenetrable obstacle.  
For more information about the previous studies see \cite{Iimpedance}, Subsections 3.2-3.3 and, in particular, 
Section 2 in which a remark on a proof in \cite{CLN} under the previous formulation of the probe method \cite{IProbe0} and \cite{IProbe2},
is given.  And it would be possible to extend Theorem 5.1 to general real-valued $V_0$ by replacing $G_y$ with a singular
solution for the Schr\"odinger equation (1.7).

\subsection{Cases covered by Theorems 1.2 and 1.3}

In this subsection, we consider an inverse obstacle problem governed by the equation
$$
\begin{array}{ll}
\displaystyle
\Delta u+k^2\,\left(a(x)+i\frac{b(x)}{k}\,\right)u=0, & x\in\Omega.
\end{array}
\tag {5.3}
$$ 
This is the case when $V_0$ and $V$ in (1.21) are given by
$$
\begin{array}{lll}
\displaystyle
V_0(x)=a_0(x)+i\frac{b_0(x)}{k}, & 
\displaystyle
V(x)=a(x)-a_0(x)+i\frac{b(x)-b_0(x)}{k}, & x\in\Omega,
\end{array}
$$
where $a$, $a_0$, $b$ and $b_0$ are real-valued and belong to $L^{\infty}(\Omega)$ ({\it no further regularity}) and $k>0$.

The equation (5.3) is coming from the time harmonic solution of the wave equation
$$
\begin{array}{ll}
\displaystyle
a(x)\partial_t^2u-\Delta u+b(x)\partial_t u=0, & x\in\Omega.
\end{array}
$$
See also \cite{CK4} for more information about the background.

We see that Assumptions 3 and 4 for this example are equivalent to the following ones.

\noindent
{\bf Assumption 3'.}  Given $F\in L^2(\Omega)$ there exists a unique weak solution $p\in H^1(\Omega)$ of
$$\left\{
\begin{array}{ll}
\displaystyle
\Delta p+k^2\,\left(a(x)+i\frac{b(x)}{k}\,\right)p=F, & x\in\Omega,\\
\\
\displaystyle
p(x)=0, & x\in\partial\Omega
\end{array}
\right.
$$
and that the unique solution satisfies
$$\displaystyle
\Vert p\Vert_{L^2(\Omega)}\le C\Vert F\Vert_{L^2(\Omega)},
$$
where $C$ is a positive constant independent of $F$.

$\quad$

\noindent
{\bf Assumption 4'.} There exists a Lebesgue measurable set $D$ of $\Bbb R^3$ having a positive measure and satisfy $D\subset\Omega$
such that $a(x)=a_0(x)$ and $b(x)=b_0(x)$ for almost all $x\in\Omega\setminus D$;

$\quad$

As a direct corollary of Theorems 1.2-1.3 we have

\proclaim{\noindent Corollary 5.1.}
Assume that Assumptions 3'-4' are satisfied.
Let $v=v(\,\cdot\,,z)$ be the complex geometrical optics solution of the equation (1.7)
given by (1.26).  Let $u$ be the unique weak solution of (5.3) with $u=v$ on $\partial\Omega$.
Then, from the indicator function (1.28) one can extract the value of the support
function of $D$ at a direction $\omega$ provided $D$ is $p$-regular with respect to $\omega$ for a $p\ge 1$
and we know one of the a-priori information (i) and (ii) below:

\noindent
(i)  $a(x)-a_0(x)$ has a positive or negative jump across $\partial D$ from $\omega$;

\noindent
(ii) $b(x)-b_0(x)$ has a positive or negative jump across $\partial D$ from $\omega$.

\endproclaim

Note that the place where $b(x)-b_0(x)\not=0$ corresponds to an unknown absorbing region
added to the background absorbing medium.
Theorem 1.3
gives us a method of estimating such a place from above.

\section{Problems remaining open}

This section describes some open issues related to enclosure method only.

\subsection{Impenetrable obstacle embedded in an absorbing medium}

Develop the enclosure method for an impenetrable obstacle embedded in an absorbing medium, that is,
the governing equation is given by the stationary Schr\"odinger equation
with a {\it fully complex-valued} background potential.  

More precisely, let $v$ and $v^*$ be given by (1.26) and (1.29), respectively.

$\quad$

{\bf\noindent Open Problem 1.}
Remove Assumption 2 and use the same form (1.28) as the indicator function for impenetrable obstacle.
Clarify, as $\tau\rightarrow\infty$ the asymptotic behaviour of this
indicator function.

$\quad$

Similarly to (1.13) we have
$$\begin{array}{ll}
\displaystyle
\left<\frac{\partial v}{\partial\nu}\vert_{\partial\Omega}-\frac{\partial u}{\partial\nu}\vert_{\partial\Omega},\,v^*\vert_{\partial\Omega}
\right>
&
\displaystyle
=-\int_{\partial D}\,\lambda(x)\,ww^*\,dS
+\int_{\Omega\setminus\overline D}\,\nabla w\cdot\nabla w^*\,dx
-\int_{\Omega\setminus\overline D}\,V_0(x)ww^*\,dx
\\
\\
\displaystyle
&
\displaystyle
\,\,\,
+\int_{\partial D}\,\lambda(x)\,vv^*\,dS
+\int_D\,\nabla v\cdot\nabla v^*\,dx
-\int_D\,V_0(x) vv^*\,dx,
\end{array}
\tag {6.1}
$$
where $u^*$ is the weak solution of (1.10) for $v$ replaced with $v^*$, $w=u-v$ and $w^*=u^*-v^*$.
Further, using the weak formulation of the governing equation of $w$, we have
$$\displaystyle
\int_{\Omega\setminus\overline D}\,\nabla w\cdot\nabla w^*\,dx
-\int_{\Omega\setminus\overline D}\,V_0(x)ww^*\,dx
=\int_{\partial D}
\left(\frac{\partial v}{\partial\nu}+\lambda(x)v\,\right)w^*\,dS.
$$
Thus (6.1) has another expression:
$$\begin{array}{ll}
\displaystyle
\left<\frac{\partial v}{\partial\nu}\vert_{\partial\Omega}-\frac{\partial u}{\partial\nu}\vert_{\partial\Omega},\,v^*\vert_{\partial\Omega}
\right>
&
\displaystyle
=\int_{\partial D}\,\left\{\left(\frac{\partial v}{\partial\nu}+\lambda(x)v\,\right)w^*-\lambda(x) ww^*\,\right\}\,dS
\\
\\
\displaystyle
&
\displaystyle
\,\,\,+
\int_{\partial D}\,\lambda(x)vv^*\,dS+\int_D\,\nabla v\cdot\nabla v^*\,dx
-\int_D\,V_0(x) vv^*\,dx.
\end{array}
\tag {6.2}
$$
Thus the main task is to study on the asymptotic behaviour of $\nabla w\cdot\nabla w^*$ over $\Omega\setminus\overline{D}$ as $\tau\rightarrow\infty$
or/and the first term of the right-hand side on (6.2).

\subsection{Extending Sini-Yoshida's another result on a penetrable obstacle}

First one has to mention another result for a penetrable obstacle in \cite{SY}. 
The governing equation is given by
$$
\begin{array}{ll}
\displaystyle
\nabla\cdot\gamma(x)\nabla u+k^2n(x)u=0, & x\in\Omega,
\end{array}
\tag {6.3}
$$
where $k$ is a fixed positive number, $\gamma$ and $n$ take the form
$$\displaystyle
\gamma(x)=
\left\{
\begin{array}{ll}
1, & x\in\Omega\setminus D,
\\
\\
\displaystyle
1+\gamma_D(x), & x\in D,
\end{array}
\right.
$$
and
$$\displaystyle
n(x)=
\left\{
\begin{array}{ll}
1, & x\in\Omega\setminus D,
\\
\\
\displaystyle
1+n_D(x), & x\in D.
\end{array}
\right.
$$
So this is the case when $V_0(x)\equiv k^2$.  
It is assumed that both $\gamma_D$ and $n_D$ are {\it real-valued}, $\gamma_D\in L^{\infty}(D)$, $n_D\in L^{\infty}(D)$ and $\text{ess.inf}_{x\in D}\gamma_D(x)>0$;
$D$ is an open subset of $\Bbb R^3$ such that $\overline{D}\subset\Omega$ and $\partial D$ is Lipschitz;
$0$ is not a Dirichlet eigenvalue for equation (6.3).

Note that $\gamma$ has a jump across $\partial D$, however, any jump condition for $n$ across $\partial D$
is not imposed.  The $D$ is a mathematical model of a penetrable obstacle embedded in a known homogeneous background medium
and appears as the perturbation term added to the {\it leading coefficient} of the governing equation.

Given a solution $v$ of the Helmholtz equation $\Delta v+k^2v=0$ in $\Omega$, 
they developed an argument to derive the following lower estimate for a fixed $p<2$:
$$\displaystyle
-\left<\frac{\partial v}{\partial\nu}\vert_{\partial\Omega}-\gamma\frac{\partial u}{\partial\nu}\vert_{\partial\Omega},\,\overline{v}\,\right>
\ge C_1\Vert \nabla v\Vert_{L^2(D)}^2-C_2\Vert v\Vert_{W^{1,p}(D)}^2-C_3\Vert v\Vert_{L^2(D)}^2,
\tag {6.4}
$$
where $u$ is the solution of (6.3) with $u=v$ on $\partial\Omega$, $w=u-v$
and $\gamma\frac{\partial u}{\partial\nu}\vert_{\partial\Omega}$ is defined by an analogous way as $\frac{\partial v}{\partial\nu}\vert_{\partial\Omega}$.

It is a consequence of the two estimates below.

\noindent
(i) The well known lower estimate which goes back to \cite{Isize}
$$\begin{array}{ll}
\displaystyle
-\left<\frac{\partial v}{\partial\nu}\vert_{\partial\Omega}-\gamma\frac{\partial u}{\partial\nu}\vert_{\partial\Omega},\,\overline{v}\,\right>
&
\displaystyle
\ge \int_{\Omega}\frac{\gamma(x)-1}{\gamma(x)}\,\vert\nabla v\vert^2\,dx-k^2\int_{\Omega}\,(n(x)-1)\,\vert v\vert^2\,dx
\\
\\
\displaystyle
&
\displaystyle
\,\,\,
-k^2\int_{\Omega}
\,(n(x)-1)\vert w\vert^2\,dx.
\end{array}
$$

\noindent
(ii) For some $p<2$ 
$$\displaystyle
\Vert w\Vert_{L^2(\Omega)}\le C\Vert v\Vert_{W^{1,p}(D)},
\tag {6.5}
$$
where $C$ is a positive constant independent of $v$.

So their contribution is an argument for the proof of (6.5).
It is divided into two parts.

{\bf\noindent Step 1.}  Prepare two facts.

$\bullet$  The Sobolev imbedding (cf. Theorem 4.12. on page 85 in \cite{AF})
$\Vert w\Vert_{L^2(\Omega)}\le C\Vert w\Vert_{W^{1,p}(\Omega)}$ with $\frac{6}{5}\le p\le 2$.

$\bullet$  The Poincar\'e inequality $\Vert w\Vert_{L^p(\Omega)}\le C\Vert \nabla w\Vert_{L^p(\Omega)}$
for $w\in W^{1,p}_0(\Omega)$ (cf. Theorem 1.4.3.4 on page 26 in \cite{Gr}).

\noindent
From those we obtain
$$\displaystyle
\Vert w\Vert_{L^2(\Omega)}\le C\Vert\nabla w\Vert_{L^p(\Omega)}.
$$
Note that this estimate is more accurate than the Poincar\'e inequality in $H^1_0(\Omega)$:
$$\displaystyle
\Vert w\Vert_{L^2(\Omega)}\le C\Vert\nabla w\Vert_{L^2(\Omega)}.
$$
Thus the problem becomes to give a good estimate of $\Vert\nabla w\Vert_{L^p(\Omega)}$ from above
in the sense that the possible upper bound is weaker than $\Vert\nabla v\Vert_{L^2(D)}$.

{\bf\noindent Step 2.}
For the purpose they employ Theorem 1 on page 198 in \cite{M} for the $L^p$-norm 
of the gradient of the solution of the elliptic problem
$$\left\{
\begin{array}{ll}
\displaystyle
\nabla\cdot\gamma\nabla r=\nabla\cdot G+H, & x\in\Omega,\\
\\
\displaystyle
r(x)=0, & x\in\partial\Omega,
\end{array}
\right.
$$
where, $p$ is in a neighbourhood of $p=2$, both $G$ and $H$ belong to $L^p(\Omega)$.  Note that the term $k^2n(x)r$ is dropped,
however, for our purpose it is no problem.
Applying this to a decomposition of $w$, one gets $\Vert\nabla w\Vert_{L^p(\Omega)}
\le C\Vert v\Vert_{W^{1,p}(D)}$ with a $p<2$ satisfying $\frac{6}{5}\le p$.
This is a brief story of the idea for the proof of (6.5).

Once we have (6.4) (and the trivial upper bound of the left-hand side on (6.4) by $C\Vert v\Vert_{H^1(D)}^2$), 
letting $v=e^{x\cdot z}$ the same as that of Subsection 5.1, showing
$$\displaystyle
\lim_{\tau\rightarrow\infty}
\frac{\Vert\nabla v\Vert_{L^p(D)}}{\Vert \nabla v\Vert_{L^2(D)}}=0
\tag {6.6}
$$
and using $\Vert\nabla v\Vert_{L^2(D)}^2\sim\tau^2\Vert v\Vert_{L^2(D)}^2$,
they established lower and upper bounds for the indicator function defined by
$$\displaystyle
I_{\omega,\vartheta}(\tau)=-\left<\frac{\partial v}{\partial\nu}\vert_{\partial\Omega}-\gamma\frac{\partial u}{\partial\nu}\vert_{\partial\Omega},\,\overline{v}\,\right>
\tag {6.7}
$$
in terms of $\Vert\nabla v\Vert_{L^2(D)}^2$.  Then some lower and upper bounds of norm $\Vert\nabla v\Vert_{L^2(D)}^2$ enabled
them to extract the value $h_D(\omega)$ from the asymptotic behaviour of the indicator function.  This gives an extension
of Theorem 1.1 in \cite{IE} which is the case $k=0$.
However, it should be noted that, as pointed out in Section 4,  their proof of (6.6) given in Lemma 3.7 in \cite{SY} 
employs the same approach explained in Section 4 and does not work for the Lipschitz obstacle case\footnote{However, using the argument
in the proof of Lemma 3.2, one could prove the validity of (6.6) for a class of obstacle including Lipschitz one.}.

$\quad$

{\bf\noindent Open Problem 2.}
It would be interested to consider the case when $n(x)$ is given by
$$\displaystyle
n(x)=
\left\{
\begin{array}{ll}
n_0(x), & x\in\Omega\setminus D,
\\
\\
\displaystyle
n_0(x)+n_D(x), & x\in D,
\end{array}
\right.
$$
where both $n_0$ and $n_D$ are essentially bounded in $\Omega$ and {\it complex-valued}.  

$\quad$

More precisely, in this case same as Theorems 1.2 and 1.3
one has to replace the indicator function (6.7) with the new one defined by
$$\displaystyle
I_{\omega,\vartheta}(\tau)=-\left<\frac{\partial v}{\partial\nu}\vert_{\partial\Omega}-\gamma\frac{\partial u}{\partial\nu}\vert_{\partial\Omega},\,v^*\,\right>,
$$
where $v^*=v(x,\overline z)$  and $v=v(x,z)$ is given by (1.26) with $V_0(x)\equiv n_0(x)$.  
Thus the problem is to clarify the aymptotic behaviour of this new indicator function.
Clearly, we will encounter the same problem as impenetrable obstacle case.

Those two open problems are in a new situation not covered by previous studies.
Such study together with the idea in Theorems 1.2 and 1.3
shall open a {\it possibility} of the application of the enclosure method to the stationary Maxwell system for an obstacle embedded in an {\it absorbing medium}.
Note that, for a non-absorbing medium, some results realizing the enclosure method exist, see \cite{KS} and \cite{KLS},
and \cite{Z} which are based on the author's previous argument in \cite{I2}.

$$\quad$$

\centerline{{\bf Acknowledgments}}

The author was partially supported by Grant-in-Aid for
Scientific Research (C)(No. 17K05331) and (B)(No. 18H01126) of Japan  Society for
the Promotion of Science.

$$\quad$$

\vskip1cm
\noindent
e-mail address

ikehataprobe@gmail.com

\end{document}